\documentclass[12pt,a4paper]{article}
\usepackage[T2A]{fontenc}
\usepackage[cp1251]{inputenc}
\usepackage{amsmath,amssymb,amsthm,amsfonts,amscd}
\begin{document}

\title{Geometric approach to stable
homotopy groups of spheres, III.  Abelian, cyclic and quaternionic
structure for mappings with singularities }
\author{P.M. Akhmet’ev}

\sloppy \theoremstyle{plain}
\newtheorem{theorem}{Theorem}
\newtheorem*{main*}{Основная Теорема}
\newtheorem*{theorem*}{Теорема}
\newtheorem{lemma}[theorem]{Lemma}
\newtheorem{proposition}[theorem]{Proposition}
\newtheorem{corollary}[theorem]{Collorary}
\newtheorem{conjecture}[theorem]{Conjecture}
\newtheorem{problem}[theorem]{Проблема}

\theoremstyle{definition}
\newtheorem{definition}[theorem]{Definition}
\newtheorem{remark}[theorem]{Remark}
\newtheorem*{remark*}{Remark}
\newtheorem*{example*}{Example}
\newtheorem{example}[theorem]{Example}
\def\aa{\dot{a}}
\def\i{{\bf i}}
\def\j{{\bf j}}
\def\a{a}
\def\b{b}
\def\bb{{\dot{b}}}
\def\k{{\bf k}}
\def\hh{{\bf \dot{h}}}
\def\e{{\bf e}}
\def\f{{\bf f}}
\def\h{{\bf h}}
\def\dd{{\dot{d}}}
\def\Z{{\Bbb Z}}
\def\R{{\Bbb R}}
\def\RP{{\Bbb R}\!{\rm P}}
\def\N{{\Bbb N}}
\def\C{{\Bbb C}}
\def\A{{\bf A}}
\def\D{{\bf D}}
\def\Q{{\mathbf Q}}
\def\QQ{{\dot{\mathbf Q}}}
\def\F{{\bf F}}
\def\J{{\bf J}}
\def\G{{\bf G}}
\def\I{{\bf I}}
\def\II{{\dot{\bf I}}}
\def\JJ{{\dot{\bf J}}}
\def\fr{{\operatorname{fr}}}
\def\st{{\operatorname{st}}}
\def\mod{{\operatorname{mod}\,}}
\def\cyl{{\operatorname{cyl}}}
\def\dist{{\operatorname{dist}}}
\def\sf{{\operatorname{sf}}}
\def\dim{{\operatorname{dim}}}
\def\dist{\operatorname{dist}}

\def\Z{{\Bbb Z}}
\def\R{{\Bbb R}}
\def\RP{{\Bbb R}\!{\rm P}}
\def\N{{\Bbb N}}
\def\C{{\Bbb C}}
\def\A{{\bf A}}
\def\D{{\bf D}}
\def\Q{{\bf Q}}
\def\i{{\bf i}}
\def\j{{\bf j}}
\def\k{{\bf k}}
\def\E{{\bf H}}
\def\F{{\bf F}}
\def\J{{\bf J}}
\def\G{{\bf G}}
\def\I{{\bf I}}
\def\e{{\bf e}}
\def\f{{\bf f}}
\def\d{{\bf d}}
\def\H{{\bf E}}
\def\fr{{\operatorname{fr}}}
\def\st{{\operatorname{st}}}
\def\mod{{\operatorname{mod}\,}}
\def\cyl{{\operatorname{cyl}}}
\def\dist{{\operatorname{dist}}}
\def\sf{{\operatorname{sf}}}
\def\dim{{\operatorname{dim}}}
\def\dist{\operatorname{dist}}
\date{}
\maketitle

\begin{abstract}
Collection of $\rm{PL}$-mappings admitting a relative abelian,
cyclic, quaternionic, bicyclic, and quaternionic-cyclic structures
are constructed.
\end{abstract}

\section*{Introduction}

A map with a target in an Euclidean space is assumed $\rm{PL}$, if a smoothness conditions do not mentioned. A generic $\rm{PL}$-map is a $\rm{PL}$-map, such that each pair of  hyperplanes spanned by the images of corresponding pair of simplexes
are transversal.
A critical point is a point, such that the restriction of the  map on an arbitrary
neighborhood of this point is not an embedding. We do not assume extra conditions for a generic $\rm{PL}$-map in  critical points.

Let us consider the groups  $\Z/2^{[s]}$, this group was defined in the
introduction of [A2] as a subgroup of the group $\Z/2 \int
\Sigma(2^{s-1})$, and the corresponding
cobordism groups of immersions (see [A2, Diagram (21)]). In [A2, Diagram  (20)] subgroups $\I_b
\times \II_b$, $\H_{b \times \bb}$, $\J_a \times \JJ_a$, $\Q
\times \Z/4$ of the groups $\Z/2^{[s]}$, $2 \le s \le 5$, are
defined and the following definitions were considered: abelian
structure (Definition 5), $\H_{b \times \bb}$--structure
(Definition 14), $\J_a \times \JJ_a$--structure (bicyclic structure)
(Definition
16), and quaternionic--cyclic structure (Definition 23) for
corresponding framed immersions. These notions are used in
Theorems 8, Lemmas 15 and 17, Theorem 25 to prove the Main Theorem in section 5.

The definitions of abelian, $\H_{b \times \bb}$--structure, $\J_a \times \JJ_a$--structure, and
quaternionic--cyclic structure of $\Z/2^{[s-1]}$--framed immersions, $s \ge 2$,
are introduced to weaken the condition of a reduction of
classifying mappings of the self-intersection $\Z/2^{[s]}$--framed immersions of the considered framed immersion, see [A2, Definitions  4, 13,
13, 22] correspondingly. Analogously, for the notion of 
quaternionic reduction see Definitions 19 in [A1]. In
the present part of the paper these notions were not considered,
the analogous relative notions were considered, and I will recall
them.

The definitions of abelian, cyclic, and quaternionic structure of framed
immersions admit  relative analogs for formal $\rm{PL}$--mappings with
singularities of the standard projective (see [A2, Definition 10]), standard $\Z/4$--lens
(see [A1, Definition 25]). The definitions of
$\H_{b \times \bb}$--structure and $\J_a \times \JJ_a$--structure of framed immersions
also admit  relative analogs for formal $\rm{PL}$--mappings with
singularities of the standard skeleton of the corresponding Eilenberg-Mac Lane spaces (see [A2, Definition 29, 31]).
The definitions of quaternionic--cyclic structure also admit
relative analogs, this analogous definition is formulated for
$\rm{PL}$--mappings with singularities of the standard skeleton
of the corresponding Eilenberg-Mac Lane space 
(see [A2,Definition 36]).

The existence of (a relative) abelian structure is formulated in
Lemma  7 of [A2], for convenience this lemma is reformulated below as
Lemma $\ref{7}$. (In the statement of this lemma below we re-denote the integer $k'$ by $k$.)

\begin{lemma}\label{7}

 For the dimensional restrictions 
\begin{eqnarray}\label{dimdim}
n-k \equiv -1 \pmod{4}, \quad k
\ge 4, \quad n \equiv 0 \pmod{2}
\end{eqnarray}
there exists a formal (equivariant) mapping $d^{(2)}: \RP^{n-k} \times \RP^{n-k} \to \R^n \times \R^n$, which  admits an abelian structure (in the sense of
[A2, Definition 10]).
 \end{lemma}


The existence of (a relative) cyclic structure is formulated in
Lemma  32 of [A1], this lemma is reformulated below as Lemma
$\ref{osnlemma1}$.

\begin{lemma}\label{osnlemma1}

A. For the dimensional restrictions 
\begin{eqnarray}\label{dim1}
n-k \equiv 1 \pmod{2},
n-3k - 10 >0, \quad n \equiv 0 \pmod {2}
\end{eqnarray} 
there exists a generic $\rm{PL}$-mapping
$d: \RP^{n-k} \to \R^n$ (with singularities) with a marked closed component of the self-intersection, 
for which the formal extension
$d^{(2)}$ admits a cyclic structure (in the sense of Definition  [A1, Definition 24]). 

B. For the dimensional restrictions 
\begin{eqnarray}\label{dimdimdim}
k \ge 5, \quad n - k \equiv 0\pmod{4} 
\end{eqnarray}
there exists a formal mapping
$d^{(2)}$  with formal self-intersection along a marked closed component $N_{\a}$, which admits a cyclic structure (in the sense of 
[A1, Definition 24]).
\end{lemma}

\begin{remark*}
Lemma $\ref{osnlemma1}$ for the proof of the main result of [A1] is not used.
\end{remark*}

The existence of (a relative) quaternionic structure is claimed in
[A1, Lemma  33] and is reformulated below. this lemma is reformulated below as Lemma
$\ref{osnlemma2}$ (in this lemma we re-denote the mapping $c$ by $d_1$).

\begin{lemma}\label{osnlemma2}
For $n=4k+(2^{\sigma}-1)$, $n=2^{\ell}-1$, $\ell \ge 7$, $\sigma =
\left[ \frac{\ell-1}{2}\right]$,
 then there
exists a generic $\rm{PL}$-mapping $d_1: S^{n-2k}/\i \to \R^n$
with singularities admitting a quaternionic structure in the sense
of [A1, Definition  $25$].
\end{lemma}

The existence of a relative 
$\H_{b \times \bb}$--structure in the sense of [A2, Definition $29$]
is formulated in [A2, Proposition 30].

The existence of a relative 
$\J_a \times \JJ_a$--structure in the sense of [A2, Definition $31$]
is formulated in [A2, Proposition 32].

The existence of a relative 
$\Q \times \Z/4$--structure in the sense of [A2, Definition $36$]
is formulated in [A2, Lemma 37].
%

In this part of the paper we shall prove all the results
formulated above from a unified point of view. The possibility of
such an approach in the case of cyclic structure was discovered by
Prof. A.V.Chernavsky at the end of the last century, and by Dr.
S.A.Melikhov (2005) in the case of quaternionic structure.
Preliminary results for cyclic and $\H_{b \times \bb}$--structure
in the case of weaker restrictions on the codimension of the
immersion, are given in the papers  [Akh1], [Akh2].
\[  \]

Let us formulate a number of remarks, which seem to be of
interest.

1. It is not, in general, possible to formulate the notion of
abelian structure (and analogous notions considered above) in
terms of the  reduction of a classifying mapping to the classifying
subspace of a corresponding abelian subgroup. For example, in the case
$n=62$ there is, as proved in [M], an obstruction to the reduction
of the classifying mapping for the self-intersection manifold of
an immersion $f: M^{n-1} \looparrowright \R^n$ into  classifying
subspace $K(\I_b \times \II_b,1) \subset K(\Z/2^{[2]},1)$ of the
abelian subgroup.

2. For the construction of cyclic and quaternionic structure for
immersions (relative cyclic and quaternionic structures for
$\rm{PL}$--mappings with singularities) only double
self-intersection points of immersions (of $\rm{PL}$--mappings) are
considered. Alternatively, in the paper [E] (this paper, as was
noted in [A1],[A2], is the foundation of our construction)
self-intersection points of an arbitrary multiplicity were
considered. In particular, it is interesting to define and to
study a quaternionic structure, related with quadruple points
manifolds of skew-framed immersions.

3. The construction of quaternionic structure in Lemma
$\ref{osnlemma2}$ does not require the Massey  embedding
$S^3/\Q \subset \R^4$ [M], see also [Me].  Such an embedding was known earlier
to W.Hantzsche [He]. By means of such
an approach, it might be possible to weaken the dimensional
restrictions in Lemma $\ref{osnlemma2}$.
For example, the Massey embedding allows to generalize Lemma
$\ref{osnlemma2}$ for maps in the range $\frac{4}{5}$ (for maps $M^{m} \to \R^n$, $\frac{m}{n} \le \frac{4}{5}$).
This means
that one may consider an extra two quadratic extensions of the
quaternionic group as the structure group of framing of
immersions.

Note that in [A1]  the cases $n=15$, $n=31$ and $n=63$
were not considered. Additional arguments, in particular, might yield a proof of
the last cases in the Adams Theorem on Hopf invariants,
and clarify the remaining case in dimension $126$ not covered by the
Hill-Hopkins-Ravenel Theorem on Kervaire invariants.

\section{Auxiliary mappings}

Cтроятся
вспомогательные отображения. В Лемме $\ref{7}$ вспомогательное
отображение $c_0$ для отображение $d_0$; в Лемме $\ref{osnlemma1}$
вспомогательные отображения $\hat c$, $c$ для отображения
$d$; в Лемме  $\ref{osnlemma2}$ вспомогательные отображения
$c_1$, $\tilde c_1$ для отображения $d_1$.

 We start by
construction of auxiliary mappings. In Lemma $\ref{7}$ this is
axillary mapping $c_0$ for the mapping
 $d_0$; in Lemma $\ref{osnlemma1}$
there are axillary mappings $\hat c$, $c$ for the mapping
 $d$; in Lemma $\ref{osnlemma2}$
there are axillary mappings $c_1$, $\tilde c_1$  for the mapping
 $d_1$.

The transformation in Lemma
$\ref{osnlemma1}$ to the required formal (equivariant) mapping $d^{(2)}$ from the mapping $c$
 is given by an approximation, which is constructed in Lemma  $\ref{Ycirc}$.

To proof the mentioned lemmas and propositions we introduce on the
singular set of auxiliary mappings the coordinate system called
\emph{angle-momentum}. By means of this coordinate system in
Lemmas $\ref{lemma280}$,$\ref{lemma291}$.
The configuration space in Lemma $\ref{lemma280}$  is defined as finite-dimensional  resolution
spaces for the singularity of the mapping $c$. In
Lemma $\ref{lemma291}$ the resolution
spaces is much simpler, because the mapping under investigation is
close to stable.

\subsubsection*{Construction of an axillary mapping  $c_0: \RP^{n-k} \to
\R^n$ in Lemma  $\ref{7}$}

Denote by $J_0$ the standard $(n-k)$--dimensional sphere of
codimension $k$ in $\R^n$, which is represented as the join of
\begin{eqnarray}\label{r}
\frac{n-k+1}{2}=r 
\end{eqnarray}
copies of the circle $S^1$. We denote the
standard embedding of $J_0$ into $\R^n$ by
\begin{eqnarray}\label{iJ}
i_{J_0}: J_0 \subset
\R^n.
\end{eqnarray}

A mapping  $p'_0: S^{n-k} \to J$ is obtained as a result of taking
the join of $r$ copies of the standard double covering $S^1 \to
\RP^1$. The standard antipodal action $\I_d \times S^{n-k} \to
S^{n-k}$ (here and below for notations of the group $\I_d$ etc.
see the first part of the section 2 in [A1]) commutes with the mapping
$p_0$. Hence, there results a mapping with ramification $p'_0:
\RP^{n-k} \to J_0$.
The required mapping
$c_0: \RP^{n-k} \to \R^n$  is defined by means of the following composition:
$i_{J_0} \circ p_0$.

\subsubsection*{Construction of axillary mappings  $c: \RP^{n-k} \to
\R^n$, $\hat c: S^{n-k}/\i \to \R^n$ in Lemma $\ref{osnlemma1}$}

The mapping $ p': S^{n-k} \to J $ is well defined as the join of
$r$ (see $(\ref{r})$) copies of the standard 4-sheeted coverings $S^1 \to S^1 / \i$. The
standard action $\I_a \times S^{n-k} \to S^{n-k}$ commutes with
the mapping $p'$. Thus, the map $ \hat p: S^{n-k} /\i \to J$ is
well defined and the map $ p: \RP^{n-k} \to J$ is well defined
as the composition $\hat p \circ \pi: \RP^{n-k} \to J$ of the standard double covering $ \pi:
\RP^{n-k} \to S^{n-k} / \i$ with the map $ \hat p$.

The required mapping $c$ is defined by the formula 
\begin{eqnarray}\label{c}
i_{J} \circ p : \RP^{n-k} \to J \subset \R^n. 
\end{eqnarray}
The required mapping $\hat c$ is defined by the formula 
\begin{eqnarray}\label{hatc}
i_{J}
\circ \hat p : S^{n-k}/\i \to \R^n.
\end{eqnarray}

\subsubsection*{Construction of axillary mappings $c_1: S^{n-2k+2^{\sigma-1}}/\i \to \R^n$,
$\tilde c_1: S^{n-2k}/\i \to \R^n$ }

Let a positive integer parameter $k$ and a positive integer $n$
are given as in Lemma $\ref{osnlemma2}$.
Let us denote by $J_1$ a $(n-2k+2^{\sigma-1})$--dimensional
polyhedron (the equation  $n-2k+2^{\sigma-1} =
\frac{n-1}{2} + 2^{\sigma}$ is satisfied), this polyhedron is defined as the
join of
\begin{eqnarray}\label{r_1}
\frac{n+1}{2^{\sigma+1}}+1=r_1
\end{eqnarray}
copies of the standard quaternionic lens space
$S^{2^{\sigma}-1}/\Q$. Below we shall used the following notation
$n_{\sigma} = 2^{\sigma}-1$, as in  [A1] and $m_{\sigma} =
2^{\sigma}-2$, as in  [A2]). By the Hirsch Theorem an embedding
 $i_{\Q}: S^{n_{\sigma}}/\Q \subset \R^{n_{\sigma}-3}$ is well defined.

Assuming $n=4k+2^{\sigma}-1$, $\ell \ge 7$ the embedding $J_1
\subset \R^n$, as the join of $r_1$ copies of the embedding
$i_{\Q}$, is well defined; let us denote this embedding by
$i_{J_1}: J_1 \subset \R^n$ (comp. with the mapping in [Lemma 35, A2].

The mapping $p'_1: S^{n-2k+n_{\sigma-1}-1} \to J_1$ is well
defined as the join of $r_1$ copies of the standard coverings
$S^{n_{\sigma}} \to S^{n_{\sigma}}/\Q$. The action $\Q \times
S^{n-2k+n_{\sigma-1}-1} \to S^{n-2k+n_{\sigma-1}-1}$ is
well defined as the standard diagonal action, given by (23)-(25)
in [A1], this action commutes with the mapping $p'_1$.

 Thus, the map $\hat p_1: S^{n-2k+n_{\sigma-1}-1}/\Q \to J_1$ is well defined and the map
 \begin{eqnarray}\label{p_1}
p_1 \cong \hat p_1 \circ \pi_1: S^{n-2k+n_{\sigma-1}-1}/\i \to
J_1,
\end{eqnarray}
  as the composition of the standard double
covering $\pi_1:
S^{n-2k+n_{\sigma-1}-1}/\i \to S^{n-2k+n_{\sigma-1}-1}/\Q$ with the map $ \hat
p_1$.

Define the required mapping
 $c_1$ as the composition
$i_{J_1} \circ p_1: S^{n-2k+n_{\sigma-1}-1}/\i \to
S^{n-2k+n_{\sigma-1}-1}/\Q \to J_1 \subset \R^n$. Consider the submanifold  $i: S^{n-2k}/\i \subset
S^{n-2k+n_{\sigma-1}-1}/\i$, this submanifold is in  general position with respect to strata
of the manifold $S^{n-2k+n_{\sigma-1}-1}/\i$, the strata are determined by the join structure.
Define the mapping
\begin{eqnarray}\label{tildep_1}
\tilde p_1 \cong \hat p_1 \circ \pi_1 \circ i: S^{n-2k}/\i \subset
S^{n-2k+n_{\sigma-1}-1}/\i \to J_1.
\end{eqnarray}
Define the required mapping
 $\tilde c_1$ as the composition 
 \begin{eqnarray}\label{tildec_1}
\tilde c_1: S^{n-2k}/\i \subset S^{n-2k+n_{\sigma-1}-1}/\i
\stackrel{c_1}{\longrightarrow} \R^n.
\end{eqnarray}

\section{Configuration spaces and singularities}

\subsubsection*{Subspaces and factorspaces of the 2-configuration space for
$ \RP^{n-k} $, related with the axillary mapping $c$ in Lemma $\ref{7}$}

In [A1, Section 3 (46)] the space $ \Gamma $, its double covering $
\bar \Gamma $, and the structural mapping $ \eta_{\Gamma}:
\Gamma \to K (\D, 1) $ were defined. 
The space $ \Gamma $ is a
manifold with boundary. Denote the interior of this manifold by $
\Gamma_{\circ} $. The restriction of the structural map $
\eta_{\Gamma} $ on $ \Gamma_{\circ} $ will be denoted by $
\eta_{\Gamma_{\circ}}: \Gamma_{\circ} \to K (\D, 1) $.

Denote by $K_{\circ} \subset \Gamma_{\circ}$ the
polyhedron of double-point singularities of the map $ p:
\RP^{n-k} \to J $, this polyhedron is defined by the formula
 $ \{[(x, y)] \in \Gamma_{\circ}, p(x) = p(y),
x \ne y \} $ (see [Formula (39),A1]). This polyhedron is equipped with a structural
mapping
\begin{eqnarray}\label{strukt}
\eta_{K_{\circ}}: K_{\circ} \to K(\D,1),
\end{eqnarray}
which is induced by the restriction of the structural mapping $
\eta_{\Gamma_{ \circ}} $ (see [A1] and below) to the subspace $K_{\circ}$.

Consider the manifold, which is defined by the
compactification of the open manifold $ \Gamma_{\circ} $ by means of diagonal  component $ \Sigma_{diag}$
(the blowing up of the diagonal is not considered).
Denote the closure of $ Cl (K_{
\circ}) $ of the polyhedron  $ K_{\circ} $  in this manifold with singularities  by $
K$.
  Denote by $Q_{diag}$
the space $Cl(K_{\circ})\setminus  K_{\circ}$. Obviously,
$Q_{diag}\subset K$. Let us call this subspace the boundary of the polyhedron $K$.

The restriction of the structure mapping
 $\eta_{K_{\circ}}$ on a regular deleted neighborhood 
$UQ_{diag\circ}$ is given by the composition of the mapping 
$\eta_{UQ_{diag\circ}}: UQ_{diag\circ} \to K(\I_b,1)$ and the mapping
$i_{\I_b,\D}: K(\I_b,1) \to K(\D,1)$. Homotopy classes of the mappings $\eta_{diag}$ and $\eta_{UQ_{diag\circ}}$
are related by the equation:
$$\eta_{diag} \circ proj_{diag} = p_{\I_b,\I_d} \circ \eta_{UQ_{diag\circ}}.$$

Note that the structural mapping of $ \eta_{K_{\circ}} $ does not
extended from $ K_{\circ} $ to the component $ Q_{diag} $ of the
boundary.  The mapping $\kappa_{diag}: Q_{diag} \to K(\I_d,1)$ is
well defined. Denote by
$U(Q_{diag})_{\circ} \subset K_{\circ}$
a small regular deleted neighborhood of $Q_{diag}$.

\subsubsection*{Subspaces and factorspaces of the 2-configuration space for
$ \RP^{n-k} $, related with the axillary mappings $c$, $\hat c$ in Lemma $\ref{osnlemma1}$}

The space $ \Gamma $,
the subspace $
\Gamma_{ \circ} \subset \Gamma $, its double coverings $
\bar \Gamma $, $\bar \Gamma_{\circ}$ were defined above. The structural mapping $\eta_{\Gamma_{ \circ}}: \Gamma_{ \circ} \to K (\D, 1)$
also were defined.

Denote by
\begin{eqnarray}\label{Sigmacirc}
\Sigma_{\circ} \subset \Gamma_{\circ}
\end{eqnarray} 
the polyhedron of double-points singularities of the map $ p:
\RP^{n-k} \to J $, this polyhedron is defined by the formula
 $ \{[(x, y)] \in \Gamma_{ \circ}, p(x) = p(y),
x \ne y \} $. This polyhedron is equipped with a structural mapping
$\eta_{\Sigma_{\circ}}: \Sigma_{\circ} \to K(\D,1),$
which is induced by the restriction of the structural mapping $
\eta_{\Gamma_{ \circ}} $ on the subspace $ \Sigma_{\circ} $.

The standard  free involution $ \i: \RP^{n-k} \to \RP^{n-k} $ is well
defined. This involution permutes points in each fiber of the
standard double covering $ \RP^{n-k} \to S^{n-k} / \i $. The space
$ \bar \Gamma_{ \circ} $ admits an involution (with fixed points)
\begin{eqnarray}\label{Ticirc}
 T_{\i\circ}
: \bar \Gamma_{ \circ} \to \bar \Gamma_{ \circ} ,
\end{eqnarray}
which is defined as the restriction of an involution 
$\i \times \i: \RP^{n-k} \times \RP^{n-k} $, constructed by the involution $
\i $ on each factor, on the subspace $ \bar \Gamma_{ \circ}
\subset \RP^{n-k} \times \RP^{n-k} $. On the quotient $\bar \Gamma_{\circ}/T=\Gamma_{\circ}$
of $ \Gamma_{\circ} $ by the another involution $ T $, which permutes the coordinates, the factorinvolution $
T_{\i \circ}: \Gamma_{ \circ} \to \Gamma_{ \circ} $ is
well defined.
 
Let us denote by $ \Sigma_{antidiag} \subset \Gamma_{\circ} $ a
subspace, called the antidiagonal, which is formed by all
antipodal pairs $ \{[(x, y)] \in \Gamma_{\circ}: x,y \in
\RP^{n-k}, x \ne y, \i (x) = y \} $. It is easy to verify that the antidiagonal $
\Sigma_{antidiag} \subset \Gamma_{\circ} $ is the set of fixed
points for the involution $ T_{\i\circ} $.

The subpolyhedron $ \Sigma_{\circ} \subset \Gamma_{\circ} $ of
multiple-points of the map $ p $ is represented by a union $
\Sigma_{ \circ} = \Sigma_{antidiag} \cup K_{\circ} $, where
$K_{ \circ} $ is an open  subpolyhedron contains all points of $
\Sigma_{\circ} $ outside the antidiagonal. The subpolyhedron $
K_{\circ} \subset \Gamma_{K_{ \circ}} $ is invariant under the
involution $ T_{\i\circ}$.

Define the restriction of the involution
 $T_{\i\circ} \vert_{K_{\circ}}$ by
$T_{K_{\circ}}$. 
The considered restriction is a free involution. Denote the factorspace 
 $K_{\circ}/ T_{K_{\circ}}$ by  $\hat K_{\circ}$.
 The restriction of the structure mapping
$\eta_{\Gamma_{\circ}}: \Gamma_{\circ}
\to K(\D,1)$ on  $K_{\circ}$  denote by 
 $\eta_{K_{\circ}}$.

Denote the closure of $ Cl (K_{
\circ}) $ of the polyhedron  $ K_{ \circ} $ (respectively, the
closure of the polyhedron $ Cl (\hat K_{ \circ}) $ polyhedron $
\hat K_ ( \circ) $) 
by $
K $ (respectively, by $ \hat K $).

 Denote by
$ Q_{diag}$ the space
$\partial \Gamma_{diag} \cap K$. Obviously, $ Q_{diag} \subset K $.
We shall call this subspace  the component of the boundary of the polyhedron $ K $.
Similarly, we denote by $ \hat Q_{diag} $ the component of the boundary of the polyhedron
$ \hat K $.

Note that the mapping $ \eta_{K} $ is not expendable to
boundary component $ Q_{diag} $. The mapping $ \kappa_{diag}:
Q_{diag} \to K (\I_d, 1) $ is well defined. Let us denote by $U(Q_{diag})_{\circ} \subset K_{\circ} $ a small regular deleted
neighborhood of  $Q_{diag} $. The projection $ proj_{diag}: U(Q_{diag})_{\circ} \to Q_{diag} $ of the regular deleted
neighborhood to  $Q_{diag} $. The restriction of the structural
mapping $ \eta_{K_{\circ}} $ to the neighborhood  $U(Q_{diag})_{\circ}$ is represented by a composition of the map $ \eta_{U(Q_{diag})_{\circ}}: U(Q_{diag})_{\circ} \to K(\I_b, 1) $ and the maps
$ i_{\I_b, \D}: K(\I_b, 1) \to K (\D, 1) $.
Homotopy classes of
maps $ \eta_{K} \vert_{Q_{diag}} $ and $ \eta_{U(Q_{diag})_{\circ}}$ satisfy the equation:
$$\eta_{diag} \circ proj_{diag} = p_{\I_b,\I_d} \circ \eta_{UQ_{diag\circ}}.$$

Let us investigate the polyhedron of singularities of an
axillary mapping $\hat c$. define the following commutative diagram of subgroups:
\begin{eqnarray}\label{HH}
\begin{array}{ccccccc}
&&&& \I_{b \times \bb} &&\\
&& \nearrow && \cap &&\\
\I_d& \subset &\I_a& \subset & \D & \subset & \E . \\
&& \searrow && \cap &&\\
&&&& \I_c &&\\
\end{array}
\end{eqnarray}

In this diagram, the inclusion
$\D \subset \E$ is a central quadratic extension of $\D$ by the element $\i$ (of the order 4),
for which  
$\i^2$ coincides with the generator $-1$ of the subgroup $\I_d \subset \D$.
The abelian groups $\E_a,
\E_{b \times \bb}, \E_c, \E_d$ are the subgroups in  $\E$, this groups are the quadratic extensions of the corresponding subgroups $\I_a, \I_b \times \II_b, \I_c, \I_d$  by means of the element $\i$. 
Note that the groups $\E_{b \times \bb}$ and $\H_{b \times \bb}$ (see above  [formula (84), A2]) are isomorphic.

The difference between the considered  groups $\E_{b \times \bb}$ and $\H_{b \times \bb}$ are the following: the representation of
$\E_{b \times \bb} \to \Z^{[3]}$ (see below [Example 16, A1]) and 
$\H_{b \times \bb} \to \Z/2^{[3]}$ (see [Diagram (85), A2]) are different. 
The kernel of the epimomorphism 
\begin{eqnarray}\label{epiE}
\E_{b \times \bb} \to \Z/2^{[3]} \to \Z/2, 
\end{eqnarray}
where $\Z/2^{[3]} \to \Z/2$ corresponds to the 
subgroup [(19),A2] of the index 2, contains an element $\i \in \E_d \subset \E_{b \times \bb}$ of the order 4
(comp. with Diagram $(\ref{140})$ below, in which $\E_d = \E_c \cap \E_{b \times \bb}$).
The kernel of the homomorphism $\H_{b \times \bb} \to \Z/2^{[3]} \to \Z/2$ coincides
with the subgroup $\I_{b \times \bb} \subset \E_{b \times \bb}$, which is an elementary 2-group.

The induced automorphism
$\chi^{[3]}: \Z/2^{[3]} \to \Z/2^{[3]}$ of the group $\E_{b \times \bb}$,
re-denoted by 
\begin{eqnarray}\label{hatchiE}
\hat \chi^{[2]}: \E_{b \times \bb} \to \E_{b \times \bb}, 
\end{eqnarray}
is defined by the formula $\hat \chi^{[2]}(\i) = \i$, where 
$\i \in \E_d$--is the generator.

The following natural mapping
$\eta_{\hat K_{\circ}}: \hat K_{\circ} \to K(\E,1)$, which corresponds to
the mapping of canonical 2-sheeted covering, is well-defined:

\begin{eqnarray}\label{140}
\begin{array}{ccccccc}
\bar K_{\circ} & \stackrel{\bar r}{\longrightarrow} & \tilde{K}_{\circ}& &K(\I_c,1)& \longrightarrow & K(\E_c,1)\\
\downarrow & & \downarrow & \longrightarrow & \downarrow & & \downarrow\\
K_{\circ} & \stackrel{r}{\longrightarrow} & \hat K_{\circ}& &K(\D,1)& \longrightarrow & K(\E,1).\\
\end{array}
\end{eqnarray}

Horizontal maps between the spaces of the diagrams we re-denote
for brevity by $ \bar \eta, \check \eta, \eta, \hat \eta $,
respectively.

\subsubsection*{Subspaces and factorspaces of the 2-configuration space for $S^{n-2k}/\i$,
related with the axillary mapping $c_1$}

The space $ \Gamma_1 $, its double covering $ \bar \Gamma_1 $, and
the structural map $ \eta_{\Gamma_1}: \Gamma_1 \to
K (\E, 1) $ was defined in [A1, Section 4, (62) and below]. The space $ \Gamma_1 $ is a manifold with boundary. Denote the interior of this manifold
by $ \Gamma_{1 \circ} $. The restriction of the structural map $ \eta_{\Gamma_1} $ to $ \Gamma_{1 \circ} $
will be denoted by $ \eta_{\Gamma_{1 \circ}}: \Gamma_{1 \circ} \to
K (\E,1) $.

Denote by $ \Sigma_{1 \circ} \subset \Gamma_{1 \ circ} $ the
polyhedron of double-points singularities of the map $ p: S^{n-2k}
\to J_1 $, this polyhedron is obtained by the blowing up of the
polyhedron $ \{[(x, y)] \in \Gamma_{1 \circ}, p(x) = p(y), x \ne y
\} $. This polyhedron is equipped with a structural mapping
\begin{eqnarray}\label{strukt2}
\zeta_{\Sigma_{1\circ}}: \Sigma_{1\circ} \to K(\E,1),
\end{eqnarray}
which
is induced by the restriction of the structural mapping $
\zeta_{\Gamma_{1 \circ}} $ on the subspace $ \Sigma_{1 \circ} $.
 
 The subpolyhedron $ \Sigma_{1 \circ} \subset \Gamma_{1 \circ} $ of
multiple-points of the map $ p_1 $ is represented by a union $
\Sigma_{1 \circ} = \Sigma_{antidiag} \cup K_{1 \circ} $, where
$K_{1 \circ} $ is an open  subpolyhedron, this subpolyhedron
contains all points of $ \Sigma_{1 \circ} $ outside the
antidiagonal. Let us denote the restriction of the structural
mapping $ \zeta_{\Gamma_{1 \circ}}: \Gamma_{1 \circ} \to K (\E, 1)
$ on $ \Gamma_{K_{1 \circ}} $ and on $ K_{1 \circ} $ by $
\zeta_{\Gamma_{1 \circ}} $ and by $ \zeta_{K_{1 \circ}} $
respectively.

 Denote the closure of $ Cl (K_{1
\circ}) $ of the polyhedron  $ K_{1 \circ} $ in $ \Gamma_{1} $ (respectively, the
closure of the polyhedron $ Cl (\hat K_{1 \circ}) $ polyhedron $
\hat K_{1 \circ} $ in $\hat \Gamma_{1} $) by $
K_{1} $ (respectively, by $ \hat K_{1} $).
  Denote by
$ Q_{antidiag}$ the space $\Sigma_{antidiag} \cap K_{1} $, denote by $ Q_{diag}$ the space
$\partial \Gamma_{diag} \cap K_{1} $. Obviously, $ Q_{diag} \subset K_{1} $,
$ Q_{antidiag} \subset K_{1} $. We shall call these subspaces  the components of the boundary of the polyhedron $K_{1}$.

Note that the structural mapping of $ \zeta_{K_{1 \circ}} $ is extended from $ K_{1 \circ} $ to the component
$ Q_{antidiag} $ of the boundary. Denote this extension by $ \zeta_{Q_{antidiag}}: Q_{antidiag} \to K(\E, 1) $.
The mapping $ \zeta_{Q_{antidiag}} $ is the composition
$ \zeta_{antidiag}: Q_{antidiag} \to K (\Q, 1) $ and the
inclusion $ i_{\Q, \E}: K (\Q, 1) \subset K (\E, 1) $.

Note that the mapping $ \zeta_{K_{1}} $ is not expendable to
boundary component $ Q_{diag} $. The mapping $ \zeta_{diag}:
Q_{diag} \to K (\I_a, 1) $ is well defined. Let us denote by $
U(Q_{diag})_{ \circ} \subset K_{1 \circ} $ a small regular deleted
neighborhood of  $Q_{diag} $. The projection $ proj_{diag}: U
(Q_{diag})_{ \circ} \to Q_{diag} $ of the regular deleted
neighborhood to $Q_{diag} $  to the central manifold is well defined.

The restriction of
the structural mapping $ \zeta_{K_{1 \circ}} $ to the neighborhood
$ U(Q_{diag})_{ \circ} $ is represented by a composition of the map $
\zeta_{UQ_{diag \circ}}: UQ_{diag \circ} \to K(\E_{b \times \bb}, 1) $ and
the maps $ i_{\E_{b \times \bb}, \E}: K(\E_{b \times \bb}, 1) \to K (\E, 1) $.

Homotopy
classes of maps $\zeta_{diag} $ and $\zeta_{UQ_{diag\circ}}$ are related by the equation: 
$$ \zeta_{diag} \circ proj_{diag} = p_{\E_{b \times \bb}, \I_a} \circ \zeta_{UQ_{diag\circ}}.$$
\[  \]







\section{Resolution spaces for singularities}

\subsection*{Resolution spaces for polyhedra  $K_{\circ}$  and $\hat K_{\circ}$}

We construct a space $RK_{\circ}$, which we call the
resolution space of the polyhedron  $K_{\circ}$. In [A2] the group $(\I_b \times \II_b)_{\chi^{[2]}} \Z$, equipped with the homomorphism
$\Phi^{[2]}: (\I_b \times \II_b)_{\chi^{[2]}} \Z \to \D$,
and the subgroup $\I_b \times \II_b \subset (\I_b \times \II_b)_{\chi^{(2)}} \Z$ are well defined. 

Consider the following diagrams: 

\begin{eqnarray}\label{16.2}
\begin{array}{ccc}
 RK_{\circ} & \stackrel {pr}{\longrightarrow} & K_{\circ}  \\
&&\\
  \phi \downarrow & & \\
&&\\
K((\I_b \times \II_b)_{\chi^{(2)}} \Z,1), &  & \\
\end{array}
\end{eqnarray}

\begin{eqnarray}\label{118.2}
\begin{array}{cccc}
RQ_{diag \circ} & \qquad \stackrel{pr}{\longrightarrow} \qquad & UQ_{diag \circ}\\
\phi \searrow & & \swarrow \eta_{diag \circ } \\
& K(\I_b \times \II_b,1), & 
\end{array}
\end{eqnarray}
where  $RQ_{diag \circ} = 
(pr)^{-1}(
UQ_{diag \circ})$.

\begin{lemma}\label{lemma28}
There exists the space
$RK_{\circ}$, which is included into the commutative diagram  $(\ref{16.2})$. 
The following diagram
$(\ref{118.2})$ determines the boundary conditions.
\end{lemma}

\subsubsection*{Resolution spaces for polyhedra   $\Sigma$ and $\hat K$}

Define a space 
$R\Sigma_{\circ}$, which is called the resolution space for the polyhedron  $\Sigma_{\circ}$, 
which is given by the formula  $(\ref{Sigmacirc})$.

The space
$R\Sigma_{\circ}$ contains two components, which is denoted by  $R\Sigma_{a}$, $RK_{b \times \bb\circ}$:
 \begin{eqnarray}\label{RK0}
R\Sigma_{a} \cup RK_{b \times \bb\circ} = R\Sigma_{\circ}.
\end{eqnarray}

The space
$R\Sigma_{a}$ is a closed polyhedron, for which the structured mapping
\begin{eqnarray}\label{phia}
\phi_a: R\Sigma_a \to K(\I_a,1)
\end{eqnarray}
is well-defined. The mapping
 $(\ref{phia})$ is included into the following commutative diagram:

\begin{eqnarray}\label{16.20}
\begin{array}{ccc}
\Sigma_{\circ} & \stackrel {pr}{\longleftarrow}& R\Sigma_{a}\\
&&\\
\downarrow \eta_{\circ} & & \downarrow \phi_a \\
&&\\
K(\D,1) & \supset &  K(\I_a,1).\\  
\end{array}
\end{eqnarray}

The space
$RK_{b \times \bb\circ}$  is a 2-sheeted covering space of the covering
 $Rr_{b \times \bb} : RK_{b\times \bb\circ} \to 
R\hat K_{b\times \bb\circ}$.

\begin{eqnarray}\label{16.23}
\begin{array}{ccccccc}
K_{\circ} & \stackrel {pr}{\longleftarrow}& RK_{b\times \bb\circ}&  \stackrel {Rr_{b\times \bb}}{\longrightarrow} &  R\hat K_{b\times \bb\circ}& \stackrel {p \hat r}{\longrightarrow} & \hat K_{\circ}  \\
&&&&&&\\
& & \downarrow \hat \phi_{b\times \bb} & & \downarrow \phi_{b\times \bb} &  & \\
&&&&&&\\
&& K(\I_{b \times \bb} \int_{\chi^{[2]}} \Z,1) & \subset &  K(\E_{b \times \bb} \int_{\hat \chi^{[2]}} \Z,1). &&  
\end{array}
\end{eqnarray}
The group
  $\E_{b \times \bb} \int_{\hat \chi^{[2]}}$, which is used in Diagram $(\ref{16.23})$
 is defined analogously to the group  $(\H_{b \times \bb}) \int_{\chi^{[3]}} \Z$,
 [Formula (68), A2], using the automorphism (involution)  $(\ref{hatchiE})$.

Denote
$(p \hat r)^{-1}(\hat
UQ_{diag\circ})$ by $R\hat Q_{diag\circ}$. The following inclusion 
$R\hat Q_{diag\circ} \subset R\hat K_{b \times \bb\circ} $ is well-defined.

Let us denote by $RQ_{diag\circ}$ the boundary of the corresponding 2-sheeted covering space over
 $R\hat Q_{diag\circ}$. The following diagram is well-defined.


\begin{eqnarray}\label{118.20}
\begin{array}{ccc}
R\hat Q_{diag\circ} & \qquad \stackrel{p\hat r}{\longrightarrow} \qquad & U\hat Q_{diag\circ} \\
&&\\
\hat \phi_{b \times \bb} \downarrow & &\hat \eta_{diag\circ} \downarrow \\
&&\\
K(\E_{b \times \bb} \int_{\chi^{[2]}} \Z,1) & \supset & K(\E_{b \times \bb},1). \\
\end{array}
\end{eqnarray}

To prove the main result of the section we will use the following lemma.

\begin{lemma}\label{lemma280}
There exists a space $R\Sigma_{\circ}$, which
is satisfies the equation $(\ref{RK0})$.

The component  $R\Sigma_a$ is equipped by the mapping
$(\ref{phia})$, which is included into the commutative diagram $(\ref{16.20})$.

The component
$RK_{b\times \bb\circ}$ is the total space of a regular 2-sheeted covering over the space 
 $R\hat K_{b\times \bb\circ}$ such that the commutative diagram  $(\ref{16.23})$
 is well-defined. Moreover, the commutative diagram
$(\ref{118.20})$, which determines boundary conditions, is well-defined. 
\end{lemma}

\subsubsection*{Resolution space for the polyhedron $\Sigma_1$}

We shall define a space
 $R\Sigma_{1\circ}$, which we call resolution
space of the polyhedron $\Sigma_1$.
The space $R\Sigma_{1\circ}$ contains two components, which is denoted by  $R\Sigma_{\Q}$, $RK_{\H_{b \times \bb}\circ}$, as follows: 
\begin{eqnarray}\label{RK1}
R\Sigma_{\Q} \cup RK_{\E_{b \times \bb}\circ} = R\Sigma_{1\circ}.
\end{eqnarray}

Let us consider the following diagrams:

\begin{eqnarray}\label{16.21}
\begin{array}{ccc}
R\Sigma_{\Q} \cup RK_{\E_{b \times \bb}\circ} & \stackrel {pr_1}{\longrightarrow} & \Sigma_1  \\
&&\\
  \phi_1 \downarrow & & \\
&&\\
K(\Q,1) \cup K(\E_{b \times \bb},1), &  & \\
\end{array}
\end{eqnarray}


\begin{eqnarray}\label{118.21}
\begin{array}{ccc}
RQ_{diag} & \qquad \stackrel{pr_1}{\longrightarrow} \qquad & Q_{diag} \\
\phi_1 \searrow & & \swarrow \zeta_{diag} \\
& K(\E_{b \times \bb},1), &
\end{array}
\end{eqnarray}
in which $RQ_{diag} = (pr_1)^{-1}(
Q_{diag})$.

The following lemma is analogous to Lemma 
$\ref{lemma280}$

\begin{lemma}\label{lemma291}

There exists a space  $RK_1$, which is satisfies the equation $(\ref{RK1})$,
an which is included in the commutative
diagram $ (\ref{16.21}) $. Moreover, the commutative
diagrams $(\ref{118.21})$ determines 
boundary conditions.
\end{lemma}

\section{Доказательствo Леммы $\ref{osnlemma1}$}

\section{Proof of Lemma $\ref{osnlemma1}$}

Let us recall that the polyhedron
 $J$ is  ${\rm{PL}}$ homeomorphic to the standard sphere $S^{n-k}$. 
 Consider the embedding $(\ref{iJ})$. Decomposes this embedding into the following composition of the standard embeddings:
$i_1: J \subset J \times \R^{k-5} \subset \R^{n-5}$,  $i_2: \R^{n-5} \subset \R^{n-1}$, $i_3: \R^{n-1} \subset \R^n$.

Consider the mapping $\hat c: S^{n-k}/\i \to
\R^n$, which is given by the formula $(\ref{hatc})$. Let us represents this mapping by the composition of the mapping  $\hat c'_1: S^{n-k}/\i \to J \times \R^{k-5}$, the inclusion  $i_2: J \times \R^{k-5} \subset \R^{n-1}$, and the standard inclusion  $i_3: \R^{n-1} \subset \R^n$.

 Define the mapping
 $\hat c_1: S^{n-k}/\i \to
\R^{n-5}$ as a result by a special  $C^1$--small ${\rm{PL}}$--deformation of the mapping  $\hat c'_1$.

Denote by
$U_{J,1} \subset \R^{n-5}$ the regular neighborhood of the embedded sphere  $J \subset J \times \R^{k-5} \subset \R^{n-5}$. Denote by 
 $proj_{J}: U_{J,1} \to J$ the orthogonal projection of a smallest neighborhood onto the central sphere
 $J$. The
${\rm{PL}}$--deformation $\hat c'_1 \mapsto \hat c_1$  is defined as a vertical deformation with respect to the orthogonal projection  $proj_{J}$.

Consider the mapping
$c_1= p \circ \hat c_1: \RP^{n-k} \to S^{n-k}/\i \to
\R^{n-5}$ and define a mapping $c'_1: \RP^{n-k} \to \R^{n-5}$ 
as the result of an additional $C^1$-small deformation $c_1 \mapsto c'_1$, which is vertical with respect to the projection $proj_{J}$, and which has the caliber $\varepsilon$, much smaller then the caliber 
$\hat \varepsilon$ of the deformation $\hat c'_1 \mapsto \hat c_1$.

Let us denote  the self-intersection polyhedron of the mapping
$c'_1$, and the open subplolyhedron the of regular self-intersection points of this map by  
\begin{eqnarray}\label{N'}
N'_{\circ} \subset N'.
\end{eqnarray}
By dimensional reasons, the mapping 
  $c'_1$  has no self-intersection points of the multiplicity  3 and more.
 Because 
the codimension $codim(\Sigma(c'_1))=k-5$, using the condition $(\ref{dim1})$ we get: $2 codim(N') > n-k$.

Because the deformation
$c_1 \mapsto c'_1$ is vertical, the polyhedron  $N'_{\circ}$ is a subpolyhedron in the polyhedron
$\Sigma_{\circ}$. Denote by
\begin{eqnarray}\label{N'K}
N'_{b \times \bb \circ} \subset  N'_{\circ}
\end{eqnarray}
an open polyhedron, which is defined by 
the inverse image of the subpolyhedron
 $(\ref{Ib})$ (see below) by the standard inclusion  $N'_{\circ} \subset \Sigma_{\circ}$.

Because $\varepsilon << \hat \varepsilon$, the subpolyhedron 
$(\ref{N'K})$ is equipped by the involution, which is induced from the involution $(\ref{Ticirc})$
 by the standard inclusion.  This involution is a free involution, because the polyhedron
$(\ref{N'K})$ does not intersects the antidiagonal.
  Let us denote by $\hat N'_{b \times \bb\circ}$ the quotient of the polyhedron 
$N'_{b \times \bb\circ}$ with respect to this involution.  
The associated 2-sheeted covering denote by
\begin{eqnarray}\label{hatN'K}
N'_{b \times \bb\circ} \to  \hat N'_{b \times \bb\circ}. 
\end{eqnarray}

The following commutative diagrams are well defined:
\begin{eqnarray}\label{250}
\begin{array}{ccc}
 N'_{\circ} & \supset & U(N'_{diag\circ}) \\
&&\\
  \downarrow \eta'_{\circ} & & \downarrow \eta'_{diag\circ} \\
&&\\
  K(\D,1)& \supset & K(\I_{b \times \bb},1),\\
\end{array}
\end{eqnarray}

\begin{eqnarray}\label{250.1}
\begin{array}{ccc}
\hat N'_{b \times \bb\circ} & \supset & U(\hat N'_{diag\circ}) \\
&&\\
  \downarrow \hat \eta'_{\circ} & &  \downarrow \hat \eta'_{diag\circ}\\
&&\\
  K(\E,1) &\supset &K(\E_{b \times \bb},1).\\
\end{array}
\end{eqnarray}

Below we shall define the required mapping
$d: \RP^{n-k} \to \R^n$ as the result of a special deformation
 $i_2 \circ c'_1 \mapsto d$.
 The deformation $i_2 \circ c'_1 \mapsto d$, generally speaking, is not a vertical deformation 
 with respect to the orthogonal projection 
 $proj_{J} \circ (\R^{n} \to \R^{n-5})$. Let us denote by
 $N_{\circ}$ an open polyhedron of self-intersection points of the mapping 
 $d$. The following subpolyhedra are well defined:  
 $N_{b \times \bb\circ} \subset N_{\circ}$, $\hat N_{b \times \bb\circ}$.
 Properties of the mapping $d$ is described in the following lemma.

\begin{lemma}\label{lemma30}
There exists a $C^0$--small  ${\rm{PL}}$-deformation $i_2 \circ c'_1 \mapsto d$, $d: \RP^{n-k} \to \R^{n-1}$, such that for the
polyhedron $N_{\circ}$ is decomposed into the union of two subpolyhedra:
\begin{eqnarray}\label{a&b}
N_{\circ} = N_{a} \cup N_{b \times \bb \circ},
\end{eqnarray}
where $N_{a}$ is closed. 

The restriction of the structure mapping
$\eta_{\circ}$ on the closed subpolyhedron  $N_{a}$ admits
a reduction, given by a mapping $\mu_a: N_a \to K(\I_a,1)$:

\begin{eqnarray}\label{inclS}
\eta_a = i_{a} \circ \mu_a: N_a \to K(\I_a,1) \subset K(\D,1).
\end{eqnarray}

The restriction of the structured map $\eta_{b \times \bb \circ}$ to the component $N_{b \times \bb \circ}$ 
is a 2-sheeted covering mapping over a mapping  $\hat \eta_{b \times \bb \circ}: \hat N_{b \times \bb \circ} \to K(\E,1)$. The mapping $\hat \eta_{b \times \bb \circ}$ admits a reduction by a mapping $\hat \mu_{b \times \bb \circ}: \hat N_{b \times \bb \circ} \to K(\E_{b \times \bb} \int_{\chi^{[2]}} \Z,1)$:
 \begin{eqnarray}\label{inclSb}
\hat \eta_{b \times \bb \circ} = \hat \Phi^{[2]} \circ \hat \mu_{b \times \bb \circ}: N_{b \times \bb \circ} \to K(\E_{b \times \bb} \int_{\chi^{[2]}} \Z,1) \to K(\E,1),
\end{eqnarray}
where $\hat \Phi^{[2]}: K(\E_{b \times \bb} \int_{\chi^{[2]}} \Z,1) \to K(\E,1)$
is a natural mapping (see an analogous [Diagram (85),A2]).
\end{lemma}

\subsubsection*{A sketch of the proof of Lemma $\ref{osnlemma1}$}

 The deformation  $i_2 \circ c'_1 \mapsto d$ will be defined, such that 
the polyhedron $(\ref{a&b})$ admits a resolution mapping:
$$ t_a \cup t_{b \times \bb \circ}: N_a \cup N_{b \times \bb \circ} \to RK_{a} \cup RK_{b \times \bb \circ}. $$
The following properties are well-defined: The mapping $t_a$ induces the following mapping
$\mu_a=\phi_a \circ t_a: N_a \to K(\I_a,1)$, which is the required mapping. 
The mapping $t_a$ induces the following mapping 
$\mu_{b \times \bb \circ} = \phi_{b \times \bb\circ} \circ t_{b \times \bb\circ}: N_{b \times \bb\circ} \to K(\I_{b \times \bb} \int_{\chi^{[2]}}\Z,1)$. This mapping
is a 2-sheeted mapping  over the second required mapping
$\hat \mu_{b \times \bb \circ}: \hat N_{b \times \bb\circ} \to K(\E_{b \times \bb} \int_{\chi^{[2]}}\Z,1)$.
An outline of the proof of Statement A of Lemma
$\ref{osnlemma1}$ is presented.
Statement B of Lemma
$\ref{osnlemma1}$ is proved analogously.

\section{Coordinate system angle-momentum
on the spaces of singularities and construction of the
resolution spaces}



\subsubsection*{The complex stratification of polyhedra $ J $, $ \Sigma $, 
$ \Sigma_{ \circ} $ by means of the coordinate
system  angle - momentum}

Let us order lens spaces, which form the join, by the integers
from 1 up to $ r $ and let us denote by $ J (k_1, \dots, k_s)
\subset J $ the subjoin,  formed by a selected set of circles
(one-dimensional lens spaces) $ S^1 / \i $ with indexes $ 1 \le
k_1 <\dots <k_s \le r $, $ 0 \ge s \ge r $. The stratification
above is induced from the standard stratification of the open
faces of the standard $ r $-dimensional simplex $ \delta^{r} $
under the natural projection $ J \to \delta^{r} $. The
preimages of vertexes of a simplex are the lens spaces $ J (j)
\subset J $, $ J (j) \approx S^1 / \i $, $ 1 \le j \le r $,
generating the join.

  Define the space $ J^{[s]} $ as a subspace of $ J $, obtained
by the union of all subspaces $ J (k_1, \dots, k_s) \subset J
$.

Thus, the following stratification
\begin{eqnarray}\label{stratJ}
J^{(r)} \subset \dots \subset J^{(1)} \subset J^{(0)},
\end{eqnarray}
of the space $
J $ is well-defined. For the considered stratum a number $ r-s $ of missed
coordinates to the full set of coordinates is called the deep of
the stratum. 

Let us introduce the following denotation:
\begin{eqnarray}\label{strat[J]}
J^{[i]} = J^{(i)} \setminus J^{(i+1)}.
\end{eqnarray}
 Denote the maximum open cell of the space $\hat p^{-1} (J
(k_1, \dots, k_s)) $ by $\hat U (k_1, \dots, k_s) \subset S^{n-k}/\i$.
 This open cell is called an elementary stratum of the
depth $ (r-s) $. A point at an elementary stratum $ U (k_1,
\dots, k_s) \subset S^{n-k} / \i $ is defined by a set of
coordinates $ (\check x_{k_1}, \dots, \check x_{k_s}, \lambda) $, where
$ \check x_{k_i} \in S^1 $ is a coordinate on the 1-sphere (circle),
covering lens space with the number $ k_i $, $ \lambda = (l_{k_1},\dots, l_{k_s}) $ is a barycentric coordinate
on the corresponding $ (s-1) $-dimensional simplex of the join.
Thus if the two sets of coordinates are identified  under the
transformation of the cyclic $ \I_a $-covering by means of the
generator, which is common to the entire set of coordinates, then
these sets define the same point on $ S^{n-k} / \i $. Points on
elementary stratum $ \hat U (k_1, \dots, k_s) $ belong in the
union of simplexes with vertexes belong to the lens spaces of the join
with corresponding coordinates. Each elementary strata $ \hat U
(k_1, \dots, k_s) $ is a base space of the double covering $ U
(k_1, \dots, k_s) \to \hat U (k_1, \dots, k_s) $, which is induced
from the double covering $ \RP^{n-k} \to S^{n-k} / \i $ by the
inclusion $ \hat U (k_1, \dots, k_s) \subset S^{n-k} / \i $.

The polyhedron  $\Sigma_{\circ}$ is split into the union of
open subsets (elementary strata),
these elementary strata are defined as the connected components of 
the inverse images of  elementary strata  $(\ref{strat[J]})$.
Denote these elementary strata by
\begin{eqnarray}\label{compstrat}
K^{[r-s]}(k_1, \dots, k_s), \qquad 1 \le s
\le r. 
\end{eqnarray}

Let us describe an elementary stratum $ K^{[r-s]} (k_1, \dots, k_s) $ by
means of the coordinate system. To simplify the notation let us
consider the case $ s = r $. Suppose that for a pair of points $
(x_1 $, $ x_2) $, defining a point on $ K^{[0]} (1, \dots, r) $,
 the following pair of points $ (\check x_1, \check x_2) $ on the covering space
  $ S^{n-k} $ is fixed, and the pair $ (\check x_1, \check x_2) $ is mapped to the pair  $ (x_1 $,
  $ x_2) $ by means of the projection of $ S^{n-k} \to \RP^{n-k} $.
  Accordingly to  the construction above, we denote by
  $ (\check x_{1, i}, \check x_{2, i}) $, $ i = 1, \dots, r $ a set of spherical coordinates of
  each point. Each such coordinate with the number $i$ defines a point on
  1-dimensional sphere
  (circle)
$ S^1_i $ with the same number $ i $,
  which covers the corresponding circle
$ J (i) \subset J $ of the join. Note that the pair of
coordinates with the common number determines the pair of points
in a common layer of the standard cyclic $ \I_a $-covering $ S^1
\to S^1 / \i $.

The collection of coordinates $ (\check x_{1, i}, \check x_{2, i})
$ are considered up to independent changes to the antipodal. In
addition, the points in the pair $ (x_1, x_2) $ does not admit a
natural order and the lift of the point in $ K $ to a pair of
points $ (\bar x_1, \bar x_2) $ on the sphere $ S^{n-k} $, is well
determined up to $8$ different possibilities. (The order of the
group $ \D $ is equal to $8$.)

An analogous construction holds for points on deeper elementary
strata $ K^{[r-s]} (k_1, \dots, k_s) $, $ 1 \le s \le r $.


\subsubsection*{The coordinate description of elementary strata of the polyhedra
$ K_{ \circ} \subset \Sigma_{\circ}$}

Let $x \in  K^{[r-s]} (k_1, \dots, k_s)$ be a point  on an
elementary stratum. Consider the sets of spherical coordinates
$\check x_{1,i}$ и $\check x_{2,i}$, $k_1 \le i \le k_s$ of the
point $ x $. For each $ i $ the following cases: a pair of $ i
$-th coordinates coincides; antipodal, the second coordinate is
obtained from first by the transformation by means of the
generator (or by the minus generator) of the cyclic cover.
Associate to an ordered pair of coordinates $ \check x_{1, k_i} $
and $ \check x_{2, k_i} $, $ 1 \le i   \le s $ the residue $ v_{k_i} = \check x_{1, k_i}(\check x_{2, k_i})^{-1} $
of a value $ +1 $, $ -1 $, $ + \i $ or $ - \i $, respectively.
It is easy to check that the collection of residues $\{v_{k_i}\}$ is changed by the following  transformation.
When  the collection of coordinates of a point is changed to the
antipodal collection, say, the collection of coordinates of the
point $ x_2 $ is changed to the antipodal collection, the set of
values of residues
 of the new pair $(\bar x_1, \bar x_2) $ on the spherical
covering is obtained from the initial set of residues by changing
of the signs:
$$\{(\check x_{1, k_i},\check x_{2, k_i})\} \mapsto \{(-\check x_{1, k_i},\check x_{2, k_i})\}, \quad \{v_{k_i}\} \mapsto \{-v_{k_i}\},$$ 
$$\{(\check x_{1, k_i},\check x_{2, k_i})\} \mapsto \{(\check x_{1, k_i},-\check x_{2, k_i})\}, \quad \{v_{k_i}\} \mapsto \{-v_{k_i}\}.$$ 
The residues of the renumbered pair of points change
by the inversion:
$$\{(\check x_{1, k_i},\check x_{2, k_i})\} \mapsto \{(\check x_{2, k_i},\check x_{1, k_i})\}, \quad \{v_{k_i}\} \mapsto \{\bar {v}_{k_i}\},$$ 
where  $v \mapsto \bar {v}$ means the complex conjugation.
  Obviously, the set of residues does not change,
if we choose another point on the same elementary stratum of the
space $ K_{\circ} $. 

Elementary strata of the space $ K (k_1, \dots, k_s) $, in
accordance with sets of residues, are divided into 3 types: $
\I_a, \I_{b \times \bb}, \I_d $. If among the set of residues are  only
residues $\{+\i,-\i\}$ (respectively, only residues $\{+1,-1\}$),
we shall speak  about the elementary stratum of the type $\I_a $
(respectively of the type $ \I_{b \times \bb} $). If among the residues are
residues from the both set $\{+\i,-\i\}$ and  $\{+1,-1\}$, we
shall speak about elementary stratum of the type $ \I_d $. It is
easy to verify that the restriction of the structure mapping $
\eta: K_{0 \circ} \to K (\D, 1) $ on an elementary stratum of
the type $ \I_a, \I_{b \times \bb}, \I_d $ is represented by the composition of
a map in the space $ K (\I_a, 1) $ (respectively in the space $ K
(\I_{b \times \bb}, 1) $ or $ K (\I_d, 1) $) with the map $i_a: K (\I_a, 1)
\to K (\D, 1) $ (respectively, with the map $ i_{b \times \bb}: K (\I_{b \times \bb}, 1)
\to K (\D, 1) $ or $ i_d: K (\I_d, 1) \to K (\D, 1) $). For
the first two types of strata the reduction of the structural
mapping (up to homotopy) is not well defined, but is defined only
up to a composition with the conjugation in the subgroups $ \I_a
$, $ \I_{b \times \bb} $.

The polyhedron $\Sigma_{\circ}$ contains the polyhedron  $K_{\circ}$ and $\Sigma_{\circ} \setminus K_{\circ}$ consists of antidiagonal elementary strata. For an arbitrary elementary antidiagonal stratum
$ K (k_1, \dots, k_s) $ the residue of the each angle coordinate is equal to  $+\i$.
A antidiagonal stratum is an elementary stratum of the type $\I_a$.
The polyhedron  $\Sigma$ is derived from  
 $\Sigma_{\circ}$ by the joining of all diagonal strata
 (on each diagonal strata the residue of an arbitrary angle coordinate is equal $+1$), which is in the boundary of the polyhedron. It is easy to verify that  $\Sigma \setminus \Sigma_{\circ}$ contains all elementary diagonal strata of the deep greater, or equal, then $1$.

Define the following open subpolyhedra
\begin{eqnarray}\label{Ia}
K_{a \circ} \subset K_{\circ} \subset \Sigma_{\circ},
\end{eqnarray}
\begin{eqnarray}\label{Ib}
K_{b \times \bb \circ} \subset K_{\circ} \subset \Sigma_{\circ},
\end{eqnarray}
\begin{eqnarray}\label{Id}
K_{d \circ} \subset K_{\circ} \subset \Sigma_{\circ}
\end{eqnarray}
as the unions of all elementary strata of the corresponding type.

The following polyhedron
\begin{eqnarray}\label{hatIb}
\hat K_{b \times \bb \circ} \subset \hat K_{\circ}
\end{eqnarray}
is defined as the base of 2-sheeted covering over the polyhedron $(\ref{Ib})$. 
The description of 
$(\ref{Ib})$ by means of the coordinates is obvious and is omitted.

\subsubsection*{Description of the structural map $\eta_{\circ}: \Sigma_{\circ} \to K(\D,1)$,
by means of the
coordinate system}

Let $x = [(x_1, x_2)]$ be a marked a point on $ K_{\circ} $,
on a maximal elementary stratum. Consider closed path $\lambda:
S^1 \to K_{\circ}$, with the initial and ending points in this
marked point, intersecting the singular strata of the depth 1 in a
general position in a finite set of points. Let $(\check x_1,
\check x_2)$ be the two spherical preimages of the point $ x $.
Define another pair $(\check x'_1,\check x'_2)$ of spherical
preimages of $ x $, which will be called  coordinates, obtained in
result of the natural transformation of the coordinates $(\check
x_1, \check x_2)$ along the path $ \lambda $.

At regular points of the path $ \lambda $ the family of pairs of
spherical preimages in the one-parameter family is changing
continuously, that uniquely identifies the  inverse images of the
end point of the path by the initial data. When crossing the path
with the strata of depth 1, the corresponding pair of spherical
coordinates with the number $ l $ is discontinuous. Since all the
other coordinates remain regular, the extension of regular
coordinates along the path at a critical moment time is uniquely
determined. For a given point $ x $ on elementary stratum of the
depth $0$ of the spaces $ K_{\circ} $ the choice of at least one
pair of spherical coordinates  is uniquely determines the choice
of spherical coordinates with the rest numbers. Consequently, the
continuation of the spherical coordinates along a path is uniquely
defined in a neighborhood of a singular point of the path.

The transformation of the ordered pair $(\check x_1, \check x_2)$
to the ordered pair $(\check x'_1,\check x'_2)$ defines an element
the group $ \D $. This element does not depend on the choice of
the path $ l $ in the class of equivalent paths, modulo homotopy
relation in the group $\pi_1(\Sigma_{\circ},x)$. Thus, the
homomorphism $\pi_1(\Sigma_{\circ},x) \to \D$ is well defined and the
induced map
\begin{eqnarray}\label{eta}
\eta_{\circ}: \Sigma_{\circ} \to K(\D,1)
\end{eqnarray}
coincides with structural mapping, which was determined earlier.
 It is easy to verify that the restriction
of the structural mapping $ \eta_{\circ} $ on the
connected components of a single elementary stratum $K_{
\circ}(1, \dots, r)$  is
homotopic to a map with the image in the subspeces $K(\I_a,1)$,
$K(\I_{b \times \bb},1)$, $K(\I_d,1)$,
which corresponds to the type and subtype elementary stratum.

\subsubsection*{Coordinate description of the canonical covering over an elementary stratum}

Consider an elementary stratum $ K^{[r-s]}(k_1, \dots, k_s) \subset
K^{(r-s)}_{\circ}$  of the depth  $ (r-s) $.
Denote by
\begin{eqnarray}\label{pi}
\pi: K^{[r-s]}(k_1, \dots, k_s) \to K(\Z/2,1)
\end{eqnarray}
 the classifying map, that is responsible for the permutation of a
pair of points around a closed path on this elementary stratum.
This mapping is called the {\it{classified}} mapping for the corresponding  2-sheeted covering.

The mapping $ \pi $  coincides with the
composition
$$ K^{[r-s]}(k_1, \dots, k_s) \stackrel{\eta}{\longrightarrow} K(\D,1) \stackrel{p}{\longrightarrow} K(\Z/2,1), $$
 where $K(\D,1)
\stackrel{p}{\longrightarrow} K(\Z/2,1)$ be the map of the
classifying spaces, which is induced by the epimorphism $ \D \to
\Z / 2 $ with kernel $ \I_c \subset \D $ .
The canonical 2-sheeted covering, which is associated with the mapping
$\pi$ let us denote  by
\begin{eqnarray}\label{bar}
\bar K^{[r-s]}(k_1, \dots, k_s) \to K^{[r-s]}(k_1, \dots, k_s).
\end{eqnarray}

With the mapping 
 $(\ref{pi})$ the following equivariant mapping is associated:
\begin{eqnarray}\label{eqvpi}
\bar{\pi}: \bar K^{[r-s]}(k_1, \dots, k_s) \to S^{\infty},
\end{eqnarray}
where the involution in the image is the standard antipodal involution.
This mapping is a 2-sheeted covering over the mapping  $(\ref{pi})$.

For an elementary strata of the type
 $\I_{b \times \bb}$ with the mapping  $(\ref{eqvpi})$ the following equivariant mapping is associated:
 \begin{eqnarray}\label{hateqvpi}
\tilde{\pi}: \tilde K^{[r-s]}(k_1, \dots, k_s) \to S^{\infty},
\end{eqnarray}
where the mapping
$\tilde K^{[r-s]}(k_1, \dots, k_s) \subset \tilde K(\E_{b \times \bb},1)$,
$(\ref{hateqvpi})$ is a 2-sheeted covering over the mapping  $(\ref{eqvpi})$.

\begin{lemma}\label{lemma32a}

The restriction of the map  $(\ref{eqvpi})$ to the canonical 2-sheeted covering over an elementary strata 
of an arbitrary type is homotopic to the following composition
\begin{eqnarray}\label{pi1a}
\bar{\pi}: \bar K^{[r-s]}(k_1, \dots, k_s) \to S^1 \subset S^{\infty}.
\end{eqnarray}
The restriction of the equivariant map  $(\ref{hateqvpi})$ to the canonical 2-sheeted covering over an elementary strata 
of the type $\E_{b \times \bb}$ is homotopic to the following composition
\begin{eqnarray}\label{hatpi1a}
\tilde{\pi}: \tilde K^{[r-s]}(k_1, \dots, k_s) \to S^1 \subset S^{\infty},
\end{eqnarray}
where $S^1 \subset S^{\infty}$ is the equivariant embedding of the standard 1-dimensional skeleton
of the classifying space.
\end{lemma}

\subsubsection*{Proof of Lemma $\ref{lemma32a}$}
Let us prove the lemma by means of explicit formulas for the mappings
$(\ref{pi1a})$ $(\ref{hatpi1a})$. An arbitrary point
$[(x_1,x_2)] \in \hat K^{[r-s,i]}(k_1, \dots, k_s)$, or $[(x_1,x_2)] \in K^{[r-s,i]}(k_1, \dots, k_s)$ 
is determined by the equivalence class of the collection of angle coordinates and the momentum coordinate. 
The structure mapping  $\eta_{\circ}$, $\hat \eta_{b \times \bb \circ}$ is determined by a transformation of angle coordinates. 
Let us define the mappings 
$(\ref{pi1a})$, $(\ref{hatpi1a})$ by the corresponding transformation of the  \emph{marked} pair of the angle coordinates. 
Below the prescribed pair of the angle coordinates for an elementary stratum of each arbitrary type is defined.   

Assume that a point
 $[(\hat{x}_1,\hat{x}_2)] \in \hat K^{[r-s]}(k_1, \dots, k_s)$ is belong to the stratum of the type $\E_{b \times \bb}$.
 Because the residue of the prescribed pair of the angle coordinates is well-defined, a non-ordered pair of the angle coordinates with the residue
 $-1$ it is convenient to denote by 
 $[(\check{x}_{1,-}, \check{x}_{2,-})]$, a pair of the angle coordinates with the residue
 $+1$ denote by  $[(\check{x}_{1,+}, \check{x}_{2,+})]$. 

The each coordinate 
 $\check{x}_{1,-}$,  $\check{x}_{2,-}$, $\check{x}_{1,+}$, $\check{x}_{2,+}$ determines the corresponding point on  $S^1$.
 It is not difficult to check, that  $\check{x}_{1,+} = \check{x}_{2,+}$, $\check{x}_{1,-} = -\check{x}_{2,-}$.
Therefore the mapping  $(\hat x_1,\hat x_2) \mapsto (\check{x}_{1,-}^{-1}\check{x}_{1,+},\check{x}_{2,-}^{-1}\check{x}_{2,+})$
transforms the points of an ordered pair into the antipodal points on  $S^1$. 
The changing of a pair of the angle coordinates to an equivalent pair, which keeps the order of the points of the pair,
does not change the equivariant mapping. The changing of the order of points in the pair transforms the equivariant mapping to the antipodal mapping.
 The constructed equivariant mapping is the required equivariant mapping $(\ref{pi1a})$ for the stratum
of the type $\E_{b \times \bb}$.

Assume a point
$[(x_1,x_2)] \in K^{[r-s,i]}(k_1, \dots, k_s)$ belongs to an elementary stratum of the type  $\I_a$ 
(including the case, when a stratum is antidiagonal). The mapping
$(\ref{pi1a})$ is determined by a transformation of the prescribed pair of the angle coordinates with the residue  $+\i$, 
which we denote (and the same time introduce an order of the pair) as  $(\check{x}_{1,+\i},\i\check{x}_{1,+\i})$. 
The mapping $(x_1,x_2) \mapsto (\check{x}_{1,+\i}^2,-\check{x}_{1,+\i}^2) $ transforms the points of the ordered pair into an antipodal
points on  $S^1$. This mapping is the required mapping
$(\ref{pi1a})$ for the elementary stratum   of the type  $\I_a$. 

Assume a point
$(x_1,x_2) \in K^{[r-s]}(k_1, \dots, k_s)$  belongs to an elementary stratum of the type $\I_d$. The mapping
$(\ref{pi1a})$ is determined by a transformation of the prescribed pair of the angle coordinates with the residue  $+\i$, 
which we denote by $[(\check{x}_{1,+\i},\i\check{x}_{1,+\i})]$. 
 The mapping $(x_1,x_2) \mapsto (\check{x}_{1,+\i})^2,-\check{x}_{1,+\i})^2 $ transforms the points of the ordered pair into an antipodal
points on $S^1$. This mapping is the required mapping
$(\ref{pi1a})$ for the elementary stratum  of the type $\I_d$.
Let us denote that the constructed mapping 
$(\ref{pi1a})$ on each elementary stratum of the type $\I_d$ is homotopic to the constant mapping.

Lemma $\ref{lemma32a}$ is proved.
\[  \]



\subsubsection*{Prescribed coordinate system and marked pair of the angle coordinates on an elementary stratum of the polyhedron
 $\hat K_{b \times \bb \circ}$}

Let us recall that the space
$\hat K_{\I_{b \times \bb} \circ}$ is the union of closures 
 $Cl(\hat K^{[r-s,i]}(k_1, \dots, k_s))$, $0
\le s \le r$ of elementary strata of the stratification $(\ref{compstrat})$ (closures are considered in the space  $\hat K_{\circ}$). 
The collection of coordinates is fixed by an ordering of the spherical
preimages $(\check{x_1},\check{x_2})$ of the marked point.  
On each elementary stratum $\hat \alpha$ of the type $\E_{b \times \bb}$ let us fix the prescribed coordinate system $\Omega(\hat \alpha)$
as follows. (In the case an equivalent class of the prescribed coordinate system of an elementary stratum depends no of an order of the preimages.)

Let us call a coordinate system a prescribed coordinate system if, 

--assuming the number of the angle coordinates is odd,
the product of residues is equal to $+1$; 

--assume that the number of the angle coordinate is even,
the number of residues $+1$ is greater then the number of residues $-1$, if
the the numbers of residues $+1$ and $-1$ coincide, the residue  with the smallest number is equal to $+1$.

The angle coordinate of the prescribed system with the residue
 $+1$ of the smallest number is called the marked coordinate on  $\hat K^{[r-s,i]}(k_1, \dots, k_s)$.

\subsubsection*{Prescribed coordinate system and marked pair of the angle coordinates on an elementary stratum of the polyhedron
 $K_{\I_a \circ}$}
 
Let us recall that the space
$\hat K_{\I_{a} \circ}$ is the union of closures 
 $Cl(K^{[r-s,i]}(k_1, \dots, k_s))$, $0
\le s \le r$ of elementary strata of the stratification $(\ref{compstrat})$ (closures are considered in the space  $K_{\circ}$). 
 On each elementary stratum $\alpha$ of the type $\I_a$ residues are $+\i$, or $-\i$. Let us define the prescribed coordinate system $\Omega(\alpha)$
as follows.

Let us call a coordinate system is the prescribed coordinate system if, 

--assuming the number of the angle coordinates is odd,
the product of residues is equal to $+\i$; 

--assume that the number of the angle coordinate is even,
the number of residues $+\i$ is greater then the number of residues $-\i$, if
the the numbers of residues $+\i$ and $-\i$ coincide, the residue with the smallest number is equal to $+\i$.

The angle coordinate of the prescribed system with the residue
 $+\i$ of the smallest number is called the marked coordinate on  $K^{[r-s,i]}(k_1, \dots, k_s)$.

\subsubsection*{Prescribed coordinate system and marked pair of the angle coordinates on an elementary stratum of the polyhedron
 $K_{\I_d \circ}$}
 
 On each elementary stratum $\alpha$ of the type $\I_d$ residues are $\{+\i,-\i, +1, -1\}$.$\{+\i,-\i\}$. Let us fix the prescribed coordinate system 
$\Omega(\alpha)$ as follows.

Let us call a coordinate system is the prescribed coordinate system if, 

--assuming the number of the angle coordinates with imaginary residues is odd,
the product of imaginary residues is equal to $+\i$; 

--assume that the number of the angle coordinate with imaginary residues is even,
the number of residues $+\i$ is greater then the number of residues $-\i$, if
the the numbers of residues $+\i$ and $-\i$ coincide, the imaginary residue  with the smallest number is equal to $+\i$.

The angle coordinate of the prescribed system with the residue
 $+\i$ of the smallest number is called the marked coordinate on  $K^{[r-s,i]}(k_1, \dots, k_s)$.
\[  \]

Let us recall that the space
$K_{\I_a \circ}$ is the union of closures 
 $Cl(K^{[r-s,i]}(k_1, \dots, k_s))$, $0
\le s \le r$ of elementary strata of the stratification $(\ref{compstrat})$
On each elementary stratum let us fix the coordinate system 
as follows.

Assume the number of the angle coordinates is odd.  Let us call a coordinate system is the prescribed coordinate system, if the sum of residues of
angle coordinates are equal to $+\i$. Assume that the number of the angle coordinate is even. Let us fixes the prescribed coordinate system
arbitrarily, namely, such that 
the residue of the pair of coordinates with the smallest number is equal to
$+\i$.

\subsubsection*{Admissible pair of neighbor strata}

Let $\beta$ be an elementary stratum (a connected component of the  space $K^{[r-s,i]}(k_1, \dots, k_s)$),
let $\alpha$ be an elementary stratum, $\alpha \subset Cl(\beta) \subset Cl(K(k_1, \dots, k_s))$, $\beta \ne \alpha$.
In this case we shall write  $\alpha \prec \beta$.

For an arbitrary $\beta \subset 
K^{[r-s,i]}(k_1, \dots, k_s)$ of the type $\I_{a}$  (correspondingly, of the type $\I_d$), let us consider an arbitrary
$\alpha$, $\alpha \prec \beta$ of the same type. Analogously, for an arbitrary $\hat \beta \subset \hat K^{[r-s,i]}(k_1, \dots, k_s)$
of the type $\E_{b \times \bb}$, let us consider an arbitrary
$\hat \alpha$, $\hat \alpha \prec \hat \beta$ of the same type.

Let us consider the prescribed coordinate system $\Omega(\beta)$ on $\beta$ and take the restriction of this coordinate system to 
$\alpha$.
Assume that the considered restriction system is prescribed on $\alpha$.  Then we shall call that the pair $(\alpha,\beta)$ is admissible.
In the case $\alpha$ and $\beta$ are of different types, we shall call that the pair $(\alpha,\beta)$ is admissible.

Assume that a pair $(\alpha,\beta)$ is not admissible. 
Take a point $b \in \beta \subset K(k_1, \dots, k_s)$ and a point $a \in \alpha$, which is closet to $b$ on $Cl(K(k_1, \dots, k_s))$.
The restriction of the prescribed coordinate system $\Omega(\beta) \vert_a$ is transformed to the prescribed system $\Omega(\alpha) \vert_a$
by one of the following transformation, which  is listed below for the strata of the each type.

A non-admissibility of a pair of strata $(\alpha,\beta)$ of the type $\I_a$
means that the transformation of $\Omega(\beta) \vert_a$ into $\Omega(\alpha) \vert_a$ is one of the following:
\begin{eqnarray}\label{prim1}
 (\check{x_1},\check{x_2}) \mapsto  (\check{x_2},\check{x_1}),
\end{eqnarray}
\begin{eqnarray}\label{prim-1}
 (\check{x_1},\check{x_2}) \mapsto  (-\check{x_2},-\check{x_1}),
\end{eqnarray}
\begin{eqnarray}\label{prim2}
 (\check{x_1},\check{x_2}) \mapsto  (-\check{x_1},\check{x_2}),
\end{eqnarray}
 \begin{eqnarray}\label{prim3}
 (\check{x_1},\check{x_2}) \mapsto  (\check{x_1},-\check{x_2}).
\end{eqnarray}

A non-admissibility of a pair of strata $(\alpha,\beta)$ of the type $\I_d$
means that the transformation of $\Omega(\beta) \vert_a$ into $\Omega(\alpha) \vert_a$ is one of the following:
\begin{eqnarray}\label{prim11}
 (\check{x_1},\check{x_2}) \mapsto  (\check{x_2},\check{x_1}),
\end{eqnarray}
\begin{eqnarray}\label{prim-11}
 (\check{x_1},\check{x_2}) \mapsto  (-\check{x_2},-\check{x_1}),
\end{eqnarray}
\begin{eqnarray}\label{prim21}
 (\check{x_1},\check{x_2}) \mapsto  (-\check{x_2},\check{x_1}),
\end{eqnarray}
 \begin{eqnarray}\label{prim31}
 (\check{x_1},\check{x_2}) \mapsto  (\check{x_2},-\check{x_1}).
\end{eqnarray}

A non-admissibility of a pair of strata $(\hat \alpha,\hat \beta)$ of the type $\E_{b \times \bb}$
means that the transformation of $\Omega(\hat \beta) \vert_a$ into $\Omega(\hat \alpha) \vert_a$ is one of the following:
\begin{eqnarray}\label{prim13}
 (\check{x_1},\check{x_2}) \mapsto  (-\check{x_2},\check{x_1}),
\end{eqnarray}
\begin{eqnarray}\label{prim-13}
 (\check{x_1},\check{x_2}) \mapsto  (\check{x_2},-\check{x_1}),
\end{eqnarray}
\begin{eqnarray}\label{prim23}
 (\check{x_1},\check{x_2}) \mapsto  (-\check{x_1},\check{x_2}),
\end{eqnarray}
 \begin{eqnarray}\label{prim33}
 (\check{x_1},\check{x_2}) \mapsto  (\check{x_1},-\check{x_2}),
 \end{eqnarray}
 \begin{eqnarray}\label{prim13i}
 (\check{x_1},\check{x_2}) \mapsto  (-\i\check{x_2},\i\check{x_1}),
\end{eqnarray}
\begin{eqnarray}\label{prim-13i}
 (\check{x_1},\check{x_2}) \mapsto  (\i\check{x_2},-\i\check{x_1}),
\end{eqnarray}
\begin{eqnarray}\label{prim23i}
 (\check{x_1},\check{x_2}) \mapsto  (-\i\check{x_1},\i\check{x_2}),
\end{eqnarray}
 \begin{eqnarray}\label{prim33i}
 (\check{x_1},\check{x_2}) \mapsto  (\i\check{x_1},-\i\check{x_2}).
\end{eqnarray}


\subsubsection*{The space $Y_{\circ}$}

Let $\alpha$, $\beta$ be elementary strata of $\Sigma_{\circ}$.
Assume that $\alpha \prec \beta$ and define the elementary $\varepsilon$-cone of a smallest stratum 
 $\alpha$ into  $\beta$ as an open neighborhood, which is defined as the open cone of a small height  $\varepsilon$, $\varepsilon <<1$,
 over the interior of the closure of the union of  all lower-dimensional $\varepsilon$-cones, which are 
inside  $Cl(\beta)$. The structure of an elementary $\varepsilon$-cone corresponds to the Euclidean structure in the $r$-simplex,
given by the corresponding momenta coordinates. The elementary cone of the strata $\alpha$ in $\beta$ denote by   
  $Con'(\alpha,\beta;\varepsilon) \subset \beta$.

For each non-admissible pair of strata 
$\alpha \prec \beta$ consider an elementary  $\varepsilon$--cone  $Con(\alpha,\beta;\varepsilon)$ and define:  
  
  -- the reduced $\varepsilon$--cone, which is denoted by  $Con^{\odot}(\alpha,\beta;\varepsilon)  \subset \beta \Sigma_{\circ}$;
 the up-reduced (correspondingly, the down-deduced)  
 $\varepsilon$--cone, which is denoted by  $Con^{\odot \uparrow}(\alpha,\beta;\varepsilon)  \subset \beta \subset \Sigma_{\circ}$ 
 (correspondingly, by  $Con^{\odot \downarrow}(\alpha,\beta;\varepsilon)  \subset \beta \Sigma_{\circ}$);
  
  --- the thickened reduced  $(\varepsilon,\varepsilon_1)$--cone, where 
\begin{eqnarray}\label{neqv}
\varepsilon_1 << \varepsilon << 1, 
\end{eqnarray}
which is denoted by 
 $Con^{\odot}(\alpha,\beta;\varepsilon,\varepsilon_1)  \subset \Sigma_{\circ}$; 
 the thickened up-reduced  $(\varepsilon,\varepsilon_1)$--cone
(correspondingly, the thickened down-reduced  $(\varepsilon,\varepsilon_1)$--cone), which is denoted by 
$Con^{\odot \uparrow}(\alpha,\beta;\varepsilon,\varepsilon_1)  \subset \beta \subset \Sigma_{\circ}$ 
(correspondingly, by  $Con^{\odot \downarrow}(\alpha,\beta;\varepsilon,\varepsilon_1)  \subset \Sigma_{\circ}$).

Let  $Con(\alpha_i,\beta; \varepsilon)$ be an arbitrary elementary cone, which is distinguished from  $Con(\alpha,\beta; \varepsilon)$, and
\begin{eqnarray}\label{prec}
\alpha \prec \alpha_i \prec \beta,
\end{eqnarray}
moreover, the pair
 $\alpha \prec \beta$ is non-admissible. Define
 $Con^{\odot \uparrow}(\alpha,\beta;\varepsilon)$ as the difference 
\begin{eqnarray}\label{elemcon}
Con(\alpha,\beta;\varepsilon) \setminus Cl(\cup_i Con(\alpha_i,\beta;\varepsilon)),
\end{eqnarray}
where
$\alpha_i$ satisfies the condition  $(\ref{prec})$ and the pair  $\alpha_i \prec \beta$ is admissible. 
Assume that instead of  $(\ref{prec})$ the following equation is satisfied: 
\begin{eqnarray}\label{precdown}
\alpha_i \prec \alpha \prec \beta.
\end{eqnarray}
Define $Con^{\odot \downarrow}(\alpha,\beta;\varepsilon)$ as the difference 
\begin{eqnarray}\label{elemcondown}
Con(\alpha,\beta;\varepsilon) \setminus Cl(\cup_i Con(\alpha_i,\beta;\varepsilon)),
\end{eqnarray}
where $\alpha_i$ satisfies the condition  $(\ref{precdown})$ and the pair $\alpha_i \prec \beta$ is admissible. 
Define  $Con^{\odot}(\alpha,\beta;\varepsilon)$ as the difference  
\begin{eqnarray}\label{elemconup}
Con(\alpha,\beta;\varepsilon) \setminus Cl(\cup_i Con(\alpha_i,\beta;\varepsilon)),
\end{eqnarray}
where $\alpha_i$ satisfies the condition  $(\ref{precdown})$, or the condition  $(\ref{prec})$, and the pair  $\alpha_i \prec \beta$ 
is admissible.

Denote by 
\begin{eqnarray}\label{cdotZo}
Z^{\odot}(\varepsilon)_{\circ} \subset \Sigma_{\circ} 
\end{eqnarray}
the disjoint union 
\begin{eqnarray}\label{cup}
\cup_{\alpha \prec \beta} Con^{\odot}(\alpha,\beta;\varepsilon),
\end{eqnarray}
where 
the pair $\alpha \prec \beta$ is non-admissible.


Consider the following CW-complex:
\begin{eqnarray}\label{ZSig}
Y_a = (\Sigma_{\circ} \setminus Z^{\odot}(\varepsilon)_{\circ}) \cap \Sigma_{a}   \subset \Sigma_{\circ},
\end{eqnarray}
where $Z^{\odot}(\varepsilon)_{\circ}$ is defined by the formula  $(\ref{cdotZo})$,  $\Sigma_{a\circ}$ is defined by the formula
 $(\ref{Ia})$. Consider the CW-complex:
\begin{eqnarray}\label{ZSigd}
Y_{d} = (\Sigma_{\circ} \setminus Z^{\odot}(\varepsilon)_{\circ}) \cap K_{d}   \subset \Sigma_{\circ},
\end{eqnarray}
where
 $K_{d\circ}$ is defined by the formula $(\ref{Id})$.
Consider the CW-complex: 
\begin{eqnarray}\label{ZSigb}
Y_{\b \times \bb} = (\Sigma_{\circ} \setminus Z^{\odot}(\varepsilon)_{\circ}) \cap K_{b \times \bb}   \subset \Sigma_{\circ},
\end{eqnarray}
where  $K_{b \times \bb \circ}$ is defined by the formula $(\ref{Ib})$.
It is not difficult to check, that the formulas
 $(\ref{prim13})$-$(\ref{prim33i})$ are invariant with respect to the covering  $(\ref{hatIb})$,
and that the CW-complex $(\ref{ZSigb})$ is thyself the covering space of the corresponding 2-sheeted covering, denote this covering by 
$Y_{\b \times \bb} \to \hat{Y}_{\b \times \bb}$.

Consider the mapping 
$\eta_{\circ}:  \Sigma_{\circ} \to K(\D,1)$, which is defined by the formula  $(\ref{eta})$.
Consider the restriction of this mapping to the subspace  $(\ref{ZSig})$ and denote this restriction by 
\begin{eqnarray}\label{etaZSiga}
\eta_a:  Y_a  \to K(\D,1).
\end{eqnarray}
Analogously, denote
\begin{eqnarray}\label{etaZSigd}
\eta_{d\circ}:  Y_d \to K(\D,1).
\end{eqnarray}
Analogously, denote
\begin{eqnarray}\label{etaZSigb}
\eta_{b \times \bb\circ}:  Y_{b \times \bb} \to K(\D,1),
\end{eqnarray}
\begin{eqnarray}\label{hatetaZSigb}
\hat \eta_{b \times \bb\circ}:  \hat Y_{b \times \bb} \to K(\E,1)
\end{eqnarray}
(see the diagram  $(\ref{HH})$).

\begin{lemma}\label{redZ}

--1. The mapping  $(\ref{etaZSiga})$ admits a reduction, which is given by the mapping 
\begin{eqnarray}\label{muZSiga}
\mu_{a\circ}:  Y_a \to K(\I_a,1),
\end{eqnarray}
$i_{\I_a,\D} \circ \mu_{a\circ} = \eta_{a\circ}$.

--2. The mapping  $(\ref{etaZSigd})$ admits a reduction, which is given by the mapping 
\begin{eqnarray}\label{muZSigd}
\mu_{d\circ}:  Y_d \to K(\I_d,1),
\end{eqnarray}
$i_{\I_d,\D} \circ \mu_{d\circ} = \eta_{d\circ}$.

--3. The mapping $(\ref{etaZSigb})$  admits a reduction, which is given by the mapping 
\begin{eqnarray}\label{muZSigb}
\mu_{b \times \bb\circ}:  Y_{b \times \bb} \to K(\I_{b \times \bb},1),
\end{eqnarray}
$i_{\I_{b \times \bb},\D} \circ \mu_{b \times \bb\circ} = \eta_{b \times \bb\circ}$.
The mapping  $(\ref{muZSigb})$ is a 2-sheeted covering over the mapping 
\begin{eqnarray}\label{hatmuZSigb}
\hat \mu_{b \times \bb\circ}:  \hat Y_{b \times \bb} \to K(\E_{b \times \bb},1).
\end{eqnarray}
\end{lemma}

\subsubsection*{Proof of Lemma $\ref{redZ}$}

Let us prove Statement 1, proofs of the last statements are analogous.
Define auxiliary spaces 
 $Y_a^{\uparrow}$ (correspondingly  $Y_a^{\downarrow}$) by the same formula that the space 
$(\ref{ZSig})$, except that in the formula  $(\ref{cup})$ the union is taken over all up-reduced (correspondingly, down-deruced) elementary $\varepsilon$--cones, which are defined by the formula   $(\ref{elemconup})$ (correspondingly, by the formula  $(\ref{elemcondown})$) 
instead of the formula $(\ref{elemcon})$. For each space 
$Y_a^{\uparrow}$, $Y_a^{\downarrow}$ the analogous statement is satisfied by the construction. 
Consider the triad 
\begin{eqnarray}\label{triad}
(\Sigma_a \setminus Y_a; \Sigma_a \setminus Y_a^{\uparrow},  \Sigma_a \setminus Y_a^{\downarrow}).
\end{eqnarray}

This triad is represented by $CW$-complexes (see below the formula 
$(\ref{cdotZoe})$). 
The required mapping
 $(\ref{muZSiga})$ is defined as the gluing the two mapping on  $\Sigma_a \setminus Y_a^{\uparrow}$,  $\Sigma_a \setminus Y_a^{\downarrow})$,
which are coincided on the small space of the triad $(\ref{triad})$. Lemma  $\ref{redZ}$ is proved.
\[  \]

Define the CW-complex
\begin{eqnarray}\label{cdotZ}
CZ^{\odot}(\varepsilon)_{\circ} \supset Z^{\odot}(\varepsilon)_{\circ},
\end{eqnarray}
as the cell closure of the space 
 $(\ref{cdotZo})$: in the CW-complex  $(\ref{cdotZ})$ all open strata of the subspace 
$(\ref{cdotZo})$ are replaced by the corresponding closure, except points on the diagonal, and the attaching mapping are continuously extended.     
The following mapping, which is a resolution, is well defined: 
\begin{eqnarray}\label{rez}
R: CZ^{\odot}(\varepsilon)_{\circ} \to \Sigma_{\circ}.
\end{eqnarray}
The restriction of the mapping  $R$ on the subspace $(\ref{cdotZo})$ is an embedding.

Let us complete  coordinates of points on an elementary cone with deleted subcones
 $(\ref{elemcon})$ by all other angle- and momentum- coordinates, which are degenerated on $\beta$, 
 the additional coordinates belong to the corresponding orthogonal face (auxiliary coordinates) to the subsimplex of (principal) momenta coordinates inside the standard  $r$-simplex.
Let us define the coordinates such that the auxiliary coordinates on  $\beta$ itself is trivial, and each auxiliary coordinate belong to 
the interval  $(0,\varepsilon_1)$. Denote this thickness by 
$Con^{\odot}(\alpha,\beta;\varepsilon,\varepsilon_1)$ and let us call it the reduced $(\varepsilon,\varepsilon_1)$--cone. 
The union of all reduced  $(\varepsilon,\varepsilon_1)$--cones 
\begin{eqnarray}\label{cdotZZ}
\cup_{\alpha \prec \beta} Con^{\odot}(\alpha,\beta;\varepsilon,\varepsilon_1),
\end{eqnarray}
where the pair 
 $\alpha \prec \beta$ is not exception, denote by  $Z^{\odot}_{\circ}(\varepsilon,\varepsilon_1)$.
Take  $\varepsilon_2 << \varepsilon_1$ and denote by  
\begin{eqnarray}\label{cdotZoe}
Z^{\odot}_{\circ}(\varepsilon,\varepsilon_1) \subset \Sigma_{\circ}
\end{eqnarray}
the subspace in  $\Sigma_{\circ}$, which is defined as the union of all reduced  $(\varepsilon,\varepsilon_1)$--cones $(\ref{cdotZZ})$. 
Denote by 
\begin{eqnarray}\label{cdotZe}
CZ^{\odot}_{\circ}(\varepsilon,\varepsilon_1) \supset Z^{\odot}_{\circ}(\varepsilon,\varepsilon_1)
\end{eqnarray}
the CW-complex, which is defined as the union of the space  $(\ref{cdotZoe})$.

The following resolution mapping
\begin{eqnarray}\label{rez}
R_{\varepsilon_1}: CZ^{\odot}_{\circ}(\varepsilon,\varepsilon_1) \to \Sigma_{\circ}
\end{eqnarray}
is well-defined. The restriction of the mapping
 $R_{\varepsilon_1}$ on the subspace  $(\ref{cdotZoe})$ is an embedding. 

Denote by
\begin{eqnarray}\label{cdotZoee}
Z^{\odot}_{\circ}(\varepsilon,\varepsilon_1,\varepsilon_2), \quad \varepsilon>>\varepsilon_1>>\varepsilon_2
\end{eqnarray}
the space, which is the union of all
$\varepsilon_2$--interiors of strata of the space  $(\ref{cdotZoe})$.
 Define 
$Y_{\circ}(\varepsilon,\varepsilon_1,\varepsilon_2)$ as the space $\Sigma_{\circ}$ with the deleted subpolyhedron  $Z^{\odot}_{\circ}(\varepsilon,\varepsilon_1,\varepsilon_2)$.

Define the space
 $Y_{\circ}$ by the formula: 
\begin{eqnarray}\label{Y}
Y_{\circ} = \lim_{\longrightarrow} (\varepsilon,\varepsilon_1,\varepsilon_2)  Y_{\circ}(\varepsilon,\varepsilon_1,\varepsilon_2), \quad \varepsilon,\varepsilon_1,\varepsilon_2 \to 0,
\end{eqnarray}
where the limit is taken over the inclusions 
 $Y_{\circ}(\varepsilon,\varepsilon_1,\varepsilon_2) \subset 
Y_{\circ}(\bar{\varepsilon},\bar{\varepsilon}_1,\bar{\varepsilon}_2)$, which are satisfies the condition  $\varepsilon > \bar{\varepsilon}$,
$\varepsilon_1 > \bar{\varepsilon}_1$, $\varepsilon_2 > \bar{\varepsilon}_2$ and the inequalities
$(\ref{neqv})$.

\begin{lemma}\label{Ycirc}

--1.  The limit $(\ref{Y})$ preserves the homotopy type of the spaces.

--2.  The CW-complex   $CZ^{\odot}_{\circ}(\varepsilon,\varepsilon_1)$ is a deformation retract of the subspace 
$Z^{\odot}_{\circ}(\varepsilon)$, which is defined by the formula  $(\ref{cdotZ})$.

--3.  Определено каноническое накрытие $\overline{CZ}^{\odot}_{\circ}(\varepsilon,\varepsilon_1) \to CZ^{\odot}_{\circ}(\varepsilon,\varepsilon_1)$, которое
 индуцировано эквивариантным отображением отображением $\bar F^{\odot}: \overline{CZ}^{\odot}_{\circ} \to \bar P$, где $\bar P$-- 3-мерное клеточное пространство со свободной инволюцией $T_P$.

--4.  The restriction of the canonical 2-sheeted covering,  which is defined in --3 over the closure of the subspace $СZ^{\odot}_{\circ} \cap K_{b \times \bb\circ}$ ($K_{b \times \bb\circ}$ is defined in $(\ref{Ib})$) is equipped by a free involution with the quotient  $\widetilde{CZ}^{\odot}_{\circ}(\varepsilon,\varepsilon_1) \to \widehat{CZ}^{\odot}_{\circ}(\varepsilon,\varepsilon_1)$, which is induced by
the following equivariant mapping 
 $\tilde F^{\odot}: \widetilde{СZ}^{\odot}_{\circ} \to \tilde P$, where $\tilde P$-- is a 3-dimensional cell complex with the involution 
 $T_{\tilde P}$.
 
\end{lemma}

\subsubsection*{Proof of Lemma $\ref{Ycirc}$}
Statement --1 is evident. 

Prove Statement --2. 
Consider the inclusion  $Z^{\odot}_{\circ}(\varepsilon) \subset Z^{\odot}_{\circ}(\varepsilon,\varepsilon_1)$. 
Using the induction over the deep of strata by the standard arguments we prove that the considered subspace is deformation retract. Statement 2 is proved.

Let us prove Statement 3. Denote by
$CZ^{\odot(s)}_{\circ} \subset CZ^{\odot(s)}_{\circ}$ the polyhedron, which consists of strata of the deep $s$ and greater, denote  $CZ^{\odot(s)}_{\circ} \setminus CZ^{\odot(s+1)}_{\circ}$ by $CZ^{\odot[s]}_{\circ}$. 
The polyhedron  $CZ^{\odot[s]}_{\circ}$ is a disjoint union of strata, which 
are differences of corresponding closures of reduced cones  
$(\ref{elemconup})$.

Define the following 3-dimensional polyhedron
$\bar P$, equipped with a free involution $T_P$. Consider the disjoint union of the elementary strata of the polyhedron $\Sigma_{\circ}$ 
and denote this union by  $\cup_s \Sigma^{[s]}$. Over each component   $\Sigma^{[s]}_i$ of  $\Sigma^{[s]}$ the canonical 2-sheeted covering
which is classified by a mappings into the circle is considered in Lemma
 $\ref{lemma32a}$.  Denote the equivariant classified mapping  by 
$\bar F_P: \cup_{s,i} \bar \Sigma^{[s]}_i \to S^1_{s,i}$. 

For each non-admissible pair of elementary strata 
 $\alpha, \beta \subset \Sigma^{(s)}$, $\alpha \prec \beta$ with the coverings  $[\bar \alpha], [\bar \beta]$ we associated the standard 
 3-sphere  $S^3_{\alpha,\beta}$, equipped with the standard action  $S^1 \times S^3_{\alpha,\beta} \to S^3_{\alpha,\beta}$. 
 Let us glue to the sphere  $S^3_{\alpha,\beta}$ the two cylinders $S^1_{\alpha} \times [0,1]$, $S^1_{\beta} \times [1,0]$ 
along the components of the boundaries  $S^1_{\alpha} \times \{0\}$, $S^1_{\beta} \times \{1\}$  
to the two antipodal fibers of the Hopf bundle, which is denoted by   $(S^1_{\alpha} \cup S^1_{\beta}) \subset S^3_{\alpha,\beta}$. 
Denote the result by $\bar{P}_{\alpha,\beta}$. The components of the boundary
$S^1_{\beta} \times \{1\}$, $S^1_{\alpha} \times \{0\}$ of the CW-complex $\bar{P}_{\alpha,\beta}$ corresponds to elementary strata of the space  $\Sigma_{\circ}$.
 

Consider the following CW-complex (non-connected) which is defined as the disjoint union of the CW-complexes
 $\{\bar{P}_{\alpha,\beta}\}$. Let us standardly identifies the circles 
 $S^1_{\alpha} \times \{0\} \cup S^1_{\beta} \times \{1\}$,  which corresponds to the common elementary stratum. 
The result is a $3$-dimensional CW-complex which is denoted by  $\bar P$. This is required space, this space is equipped with the standard antipodal involution which is denoted by   $T_P$. 

Define the following 1-dimensional CW-complex 
$\bar Q \subset \bar P$ (non-connected), which is invariant with respect to the involution $T_P$, this space is given by the union of circles 
$\{S^1_{\alpha}\}$, the components of this space corresponds to the elementary strata of the space  $\bar \Sigma_{\circ}$. 
The components $\bar Q$ are equipped with a natural stratification which is denoted by  $\bar Q^{[i]}$. 
The stratification is defined as deeps of  strata.

Define the space
 $CZ^{\odot [i]}$, the components of this space corresponds to differences of reduced cones in closures of elementary strata of $\Sigma^{[i]}$ of the deep 
 $i$. The following equivariant mapping
 $\bar F_P^{[i]}: \overline{CZ}^{\odot [i]} \to \bar P$ is well-defined, the image of this mapping belongs to  $\bar Q \subset \bar{P}$. 
This equivariant mapping is defined by the formula  $\bar F_P^{[i]}: \cup_{\alpha \prec \beta} Cl(\overline{Con}^{\odot}(\alpha,\beta;\varepsilon)) \to \bar P$. Below we shall write "`mapping"' instead of "`equivariant mapping"' for short.

 Proof of Statement 3 is given by the induction. Define 
 $P^{(s)}$ as the subspace in $\bar{P}$, which is the union of $\{\bar{P}_{\alpha,\beta}\}$, where the deep of each strata is not less then then $s$.  
Over the polyhedron $P^{(s)}$ the canonical 2-sheeted covering $\bar{P}^{(s)} \to P^{(s)}$ is well-defined and this covering is equipped by the free involution which will be denoted by $T_P^{(s)}$.  Let us prove that the mapping  $\bar F_P^{(s+1)}$ is extended from  $\overline{CZ}^{\odot(s+1)}$ 
to $\bar{P}^{(s+1)}$ into a mapping 
 $\bar F_P^{(s)}$ from $\overline{CZ}^{\odot(s)}$ to  $\bar{P}^{(s)}$.

Assume that the mapping
 $\bar{F}^{(s+1)}: \overline{CZ}^{\odot(s+1)} \to \bar{P}^{(s+1)}$ is well defined, moreover this mapping satifies the following condition.
Let us mark for each reduced elementary cone of the deep not less then 
 $s+1$ the standard
$r-s-1$--dimensional torus which is determined by the momentum coordinate near the vertex of the cone.  
It is required that in a neighborhood of this marked torus the mapping
$\overline{CF}^{(s+1)}$ coincides to the standard mapping into the circle, which is constructed in Lemma  $\ref{lemma32a}$, correspondingly to the type of the strata, which contains the elementary cone.


Let us construct the mapping
$\bar{F}^{(s)}: Z^{\odot(s)} \to \bar P^{(s)}$, which satisfies the analogous conditions as the mapping  $\bar{F}^{(s+1)}$. 
Consider an arbitrary elementary stratum $\beta$ of the deep  $s$ in $\Sigma^{[s]}_{\circ}$. 
The prove is given by an induction over the decrease of the deep $j$ of strata 
$\alpha_1$, where the pair $\alpha_1 \prec \beta$ is non-admissible. 
Namely, consider in $\cup_{i} Con^{\odot}(\alpha_i,\beta;\varepsilon)$ the union of all reduced cones of the deep more them $j$.
Then we continue the mapping over this union to each elementary cone, which is constructed from the stratum   $\alpha_1$ of the deep $j$. 
The key obvious observation is the following.

\subsubsection*{Observation (H)}
Consider a triple of strata 
$\alpha_1 \prec \beta$, $\alpha_2 \prec \beta$, $\alpha_2 \prec \alpha_1$, assuming that the first two pairs are non-admissible, the deep of 
$\beta$ is equal to  $s$, the deep of  $\alpha_1$ is equal to  $j$, the deep of  $\alpha_2$ is more then  $j$. Then the pair  $\alpha_2 \prec \alpha_1$ is admissible.
\[  \]

Using the denotations introduced above consider the reduced cone
 $Con^{\odot}(\alpha,\beta;\varepsilon)$, where $\alpha \prec \beta$ is non-admissible, and consider inside this cone all smallest elementary cones
$\alpha_i$, such that the pairs $\alpha_i \prec \alpha$, $\alpha_i \prec \beta$ are non-admissible. 
Recall that the deep  of $\alpha$ is equal to  $j$,  the deep  of $\beta$ is equal to $s$, $j < s$. Let us fixes  $\delta > 0$, $\delta << \varepsilon$. Consider an open domain  $\Omega(\alpha,\beta;\varepsilon,\delta)$, which is defined as the result of the 
elimination from the cone $\beta$ of all elementary  $\varepsilon - \delta$--cones of all strata $\alpha_i$ of the deep more then $j$, such that 
the pair $\alpha_i \prec \beta$ is non-admissible, and also the pair $\alpha_1 \prec \beta$ is non-admissible.

Define the mapping
$\bar{F}_{\alpha_1,\beta}: \overline{\Omega}(\alpha,\beta;\varepsilon,\delta) \to S^3_{\alpha,\beta}$, which is called the standard. 
Consider a regular equivariant 
$\frac{\delta}{4}$--neighborhood of the strata  $\bar{\alpha}$ is the subspace
$\bar{\Omega}(\alpha,\beta;\varepsilon,\delta)$ and denote this neighborhood by $\bar{W}(\alpha)$. 

Consider the difference
$\alpha \setminus \cup_i Con^{\odot}(\alpha_i,\alpha;\varepsilon)$, where the pair $\alpha_i \prec \alpha$ is admissible, and denote this difference by
 $\alpha^{\odot}$.  Because the cone $Con^{\odot}(\alpha,\beta;\varepsilon)$ in an up-reduced cone, by the Observation (H) an arbitrary cone
 $C(\alpha,\alpha_1) \subset \alpha_1$, $\alpha \prec \alpha_1 \prec \beta$, where the pair $\alpha_1 \prec \alpha$ is non-admissible,
 has no intersection with 
$Con^{\odot}(\alpha_i,\alpha;\varepsilon)$.

Define the mapping 
$\bar{F}_{\alpha_1,\beta}$ on $\bar{W}(\alpha)$, which is in the boundary of $\bar{\alpha}^{\odot}$, as the composition 
of the equivariant projection on
$\bar \alpha$ with the mapping  $\bar F_{\alpha} \bar \alpha \to S^1_{\alpha} \subset \bar P$. 
Define the mapping  $\bar{F}_{\alpha_1,\beta}$ on a part of  $\bar{W}(\alpha)$, which is in the boundary of $\bar{W}(\alpha) \subset \overline{Con}^{\odot}(\alpha,\beta;\varepsilon,\delta) \subset \bar \beta$, as the composition of the equivariant inclusion on
$\bar \beta$ with the mapping  $\bar F_{\beta} \bar \beta \to S^1_{\beta} \subset \bar P$. The mapping
 $\bar{F}_{\alpha,\beta}$ on $\overline{\Omega}(\alpha,\beta;\varepsilon,\delta) \setminus \bar{W}(\alpha)$
is defined analogously as above. 

Define the mapping
  $\bar{F}_{\alpha,\beta}$ on $\bar{W}(\alpha)$ by the linear approximation of the prescribed boundary conditions, which are considered as the pair of complex-valued mappings into the Whitney sum of the complex line bundles.  The standard mapping 
$\bar{F}_{\alpha,\beta}: \overline{Con}^{\odot}(\alpha,\beta;\varepsilon) \to S^3_{\alpha,\beta}$
is well-defined. The standard mapping   
 $\bar{F}_{\alpha_1,\beta}$ is continuously extended into the closure  $Cl(\overline{\Omega})(\alpha,\beta;\varepsilon,\delta)$.
Denote this extension  by  $\overline{CF}_{\alpha,\beta}: Cl(\overline{\Omega})(\alpha,\beta;\varepsilon,\delta) \to \bar P$.

It is claimed:

--1. The mapping  $\overline{CF}_{\alpha_1,\beta}$ corresponds to the mapping, which is defined on previous steps of the construction
on a deeper cone 
$\overline{Con}^{\odot}(\alpha_1,\alpha;\varepsilon)$, such a cone is included into the stratum  $\alpha$, 
moreover the pair  $\alpha_1 \prec \alpha$ is non-admissible. 

--2. The restriction of the mapping $\bar{CF}_{\alpha,\beta}$  on the domain  $\overline{\Omega}(\alpha_1,\beta;\varepsilon,\delta)$ 
inside each deeper cone is agree with the mapping  $\overline{CF}_{\alpha_1,\beta}$, where $\alpha_1 \prec \alpha \prec \beta$.

Prove --1, using Observation (H). Because the pair 
$\alpha_1 \prec \alpha$ is non-admissible, the elementary cone 
$\overline{Con}(\alpha_1,\beta;\varepsilon)$ has no intersection with  $\overline{\Omega}$. The boundary condition 
over
$\alpha^{\odot}$ of the mapping 
$\overline{CF}_{\alpha,\beta}$ proves the Statement 1. 

Prove --2, using Observation (H). By the construction  the mapping
 $\overline{CF}_{\alpha_1,\beta}$ is induced by the mapping  $\bar F_{\beta}$ everywhere on
$\overline{\Omega}(\alpha,\beta;\varepsilon,\delta) \cup \overline{Con}(\alpha_1,\beta;\varepsilon-\frac{\delta}{2})$.
The mapping  $\bar F_{\beta}$ is induced by the same mapping on the considered intersection. Statement 2 is proved. 

Statement 3 is proved. Statement 4 is evident. Lemma 
$\ref{Ycirc}$ is proved.

\subsubsection*{The canonical covering over $K_{d\circ} \subset \Sigma_{\circ}$}

Consider the subspace $K_{d\circ} \subset \Sigma_{\circ}$, which is defined by the formula  $(\ref{Id})$.
The following lemma precises Lemma 
 $\ref{Ycirc}$, Statement 3. 

\begin{lemma}\label{Ybbcirc}
  
The canonical covering over the subspace   
 $K_{d\circ} \subset \Sigma_{\circ}$ 
 is induced by an equivariant mapping
 $\bar F_{d\circ}^{\odot}: K_{d\circ} \to \bar P_{d\circ}$, where $\bar P_{d\circ}$ is a 4--dimensional CW-complex, equipped with a free involution $T_{P_{d\circ}}$.
 
\end{lemma}

\subsubsection*{Proof of Lemma  $\ref{Ybbcirc}$}

Consider the subspace 
 $Y_{d\circ} \subset K_{d\circ}$,  which is defined by the formula  $(\ref{ZSig})$. 
 The canonical covering over this subspace is trivial (see the formula $(\ref{muZSigd})$).
By Lemma  $\ref{Ycirc}$, Statement 3, the canonical covering over the subspace $K_{d\circ} \setminus Y_{d\circ}$ 
is classified by a mapping into 3-dimensional CW-complex. Lemma  $\ref{Ybbcirc}$ is proved. 
\[  \]



\subsubsection*{Definition of spaces  $R\Sigma_a$, $R\hat K_{b\times \bb\circ}$ in Lemma  $\ref{lemma280}$}
Define the subspace  
\begin{eqnarray}\label{RSig}
R\Sigma_a \subset Y_{\circ},
\end{eqnarray}
which consists of strata of the type 
$\I_a$ (c. with  $(\ref{ZSig})$).
 
Define the space 
\begin{eqnarray}\label{hatRSig}
RK_{b\times \bb\circ} \subset  Y_{\circ},
\end{eqnarray}
which consists of strata of the type
$\I_{b \times \bb}$ (c. with  $(\ref{ZSigb})$).
The space 
 $(\ref{hatRSig})$ is a 2-sheeted covering space, denote the base of the covering by  $R\hat K_{b \times \bb}$.

Definitions of the mappings, which are included into the diagram 
 $(\ref{16.23})$, in particular, the mappings  $pr$, $p \hat r$, are evident.


\subsubsection*{Resolution mapping  $\phi_a: R\Sigma_a \to K(\I_a,1)$
and Proof of Lemma $\ref{lemma280}$}

Consider the restriction 
\begin{eqnarray}\label{etaaRK}
\eta_{\circ} \vert_{Y_a}: Y_a \to K(\D,1),
\end{eqnarray} 
(recall that  $R\Sigma_a=Y_a$)
of the structured mapping 
  to the
subpolyhedron $(\ref{RSig})$.
By the construction of the reduction mapping 
\begin{eqnarray}\label{etaaaRK}
\phi_{Y_a}: Y_{a} \to K(\I_a,1),
\end{eqnarray} 
of the mapping $(\ref{etaaRK})$
is well defined: 
 $\eta_{\circ} \vert_{Y_a}= i_{\I_a} \circ \phi_{Y_a}$.
Lemma $\ref{lemma280}$ is proved. 
\[  \]

\subsubsection*{Last step of the proof of Lemma  $\ref{lemma30}$; the deformation  $i_2 \circ c'_1 \mapsto d$}


Denote the standard orthogonal projection 
$\R^{n} \to \R^{n-5}$ by  $F$. Assuming the dimensional restriction $(\ref{dim1})$, using Lemma  $\ref{lemma30}$ and
Lemma $\ref{Ycirc}$ Statement.3, let us define a vertical  lift of the mapping $F$:  $i_1 \circ i_2 \circ c'_1 \mapsto d: J \to \R^n$, $i_1: \R^{n-k-5} \subset \R^{n}$ (see denotations in Lemma  
$(\ref{lemma30})$, such that the self-intersection polyhedron $\N(d)$ is contained into the polyhedron  $Y_{\circ}$, see. $(\ref{Y})$.  
Self-intersection points of the mapping  $d$ are divided into two closed subpolyhedra correspondingly with 
the required formula
$(\ref{a&b})$. The required mappings  $\hat \mu_{b \times \bb \circ}$, $\mu_a$ are induced from the mappings, which are constructed in Lemma  $\ref{lemma280}$. Lemma $\ref{lemma30}$ is proved. 

\subsubsection*{Proof of Lemma 
 $\ref{osnlemma1}$}
 
Assuming the dimensional restriction
 $(\ref{dimdimdim})$ let us consider an axillary mapping  $(\ref{c})$ and the mapping
$F \circ c: \RP^{n-k} \to J \subset \R^{n-5}$. Consider the formal (equivariant) mappings
 $(F \circ c)^{(2)}$, $c^{(2)}$, which are defined as the formal extensions of the corresponding mappings. 
The polyhedrons of the (formal) self-intersection of the formal mappings  $(F \circ i_1 \circ с)^{(2)}$ and
$c^{(2)}$ coincide. The equivariant deformation of the formal (equivariant) mapping  $(F \circ i_1 \circ с)^{(2)}$ into the formal (equivariant)
mapping  $d^{(2)}$, which is vertical along $F^{(2)}$ is defined as in Lemma  $\ref{lemma30}$.
 

Let us prove two conditions in the statement of [Lemma 27, A1]. Condition 1 is, obviously, well proved,
namely, the restriction of the mapping  $\eta_{\circ}$
to the 
marked component $N_a$ admits a cyclic reduction, given by  $\mu_a$.

Let us prove Condition 2 in [Lemma 27, A1], which is formulated for the component
$N_{b \times \bb\circ}$.  For the convenience let us write-down this condition: 
\begin{eqnarray}\label{usll}
0=(p_{\I_c,\I_d} \circ \bar \eta)_{\ast}([\bar N_{b \times \bb}]) \in  H_{n-2k}(K(\I_d,1);\Z/2).
\end{eqnarray}

Assume that the polyhedron 
 $N_{b \times \bb\circ}$ is closed (let us remain that in this case the lower index $\circ$ in omitted) 
 and the mapping  $\eta$ admits a reduction
\begin{eqnarray}\label{redetabb}
\eta_{b \times \bb}: N_{b \times \bb} \to K(\I_{b \times \bb},1).
\end{eqnarray}
In this case the formula 
 $(\ref{usll})$ is satisfied, because the composition  
$$\bar \eta_{b \times \bb}: \bar N_{b \times \bb} \to K(\I_d,1)$$ 
is the composition of a mapping $N_{b \times \bb} \to K(\I_d,1)$
with the standard 2-sheeted covering
$$ \bar N_{b \times \bb} \to N_{b \times \bb} \to K(\I_{b \times \bb},1) \to K(\I_d,1),$$ 
where the mapping $K(\I_{b \times \bb},1) \to K(\I_d,1)$ is induced by the homomorphism 
$\I_{b \times \bb} \to \I_d$ with the kernel $\I_b \subset \I_{b \times \bb}$.

Assume that the polyhedron
 $N_{b \times \bb\circ}$ is not closed, and the mapping 
$\eta_{\circ}$ admits a reduction  $(\ref{redetabb})$ with the prescribed boundary conditions. 
The formula 
 $(\ref{usll})$ is rewritten as follows:
\begin{eqnarray}\label{usll2}
0=(p_{\I_c,\I_d} \circ \bar \eta_{b \times \bb \circ,\circ})_{\ast}([C\bar N_{b \times \bb\circ}]) \in  H_{n-2k}(K(\I_d,1);\Z/2).
\end{eqnarray}
The difference between the formulas 
 $(\ref{usll2})$ and $(\ref{usll})$ is following: if the polyhedron 
$N_{b \times \bb \circ}$ is non-closed, then the polyhedron  $\bar N_{b \times \bb\circ}$
is also non-closed. Therefore the polyhedron
$\bar N_{b \times \bb\circ}$ have to be compactified into a closed by a gluing of the cone of the canonical 2-sheeted cover 
$\bar N_{b \times \bb \circ} \to  N_{b \times \bb \circ}$ over the boundary. 
The result is a closed polyhedron, which is denoted in the formula
$(\ref{usll2})$ by $С\bar N_{b \times \bb \circ}$. The polyhedron
$С\bar N_{b \times \bb \circ}$ is the covering space of the 2-sheeted covering 
$С\bar N_{b \times \bb \circ} \to C N_{b \times \bb \circ}$, which corresponds to the subgroup 
$\I_{\bb} \subset \I_{b \times \bb}$ of the index 2. Therefore, as in the previous case, the cycle
$p_{\I_c,\I_d} \circ \bar \eta_{b \times \bb \circ}: C\bar N_{b \times \bb \circ} \to K(\I_d,1)$ is a boundary.

Let us consider a general case: the polyhedron  $N_{b \times \bb \circ}$ is non-closed and the mapping
$\eta_{\circ}$ admits a reduction 
$$\eta_{b \times \bb\circ}: N_{b \times \bb\circ} \to K(\I_{b \times \bb} \int_{\chi^{[2]}} \Z,1)$$
with prescribed boundary conditions.

By the assumption the following mapping
$$\hat \eta_{b \times \bb\circ}: \hat N_{b \times \bb\circ} \to K(\E_{b \times \bb} \int_{\chi^{[2]}} \Z,1)$$
is well-defined. Consider the 2-sheeted covering over the structure mapping, which we denote by 
$$\tilde \eta_{b \times \bb\circ}: C\widetilde{N}_{b \times \bb \circ} \to K(\E_{d} \times \Z,1).$$ 

Let us recall, that respectively to the diagram 
$(\ref{140})$, the 2-sheeted covering mapping $\tilde \eta_{b \times \bb \circ}$ over $\eta_{b \times \bb \circ}$ 
is totally defined by the subgroup of the index 2: 
\begin{eqnarray}\label{Ebb}
\E_{d} \times \Z \subset \E_{b \times \bb} \int_{\hat \chi^{[2]}} \Z.
\end{eqnarray}

The formula
$(\ref{usll2})$ is equivalent to the following condition: the homology class
\begin{eqnarray}\label{usll4}
(p_{\E_{d} \times \Z,\E_d} \circ \tilde \eta_{b \times \bb\circ})_{\ast}([C\widetilde{N}_{b \times \bb \circ}]) \in  H_{n-2k}(K(\E_d,1);\Z)
\end{eqnarray}
is even.

By the representation 
$\E_{b \times \bb} \int_{\hat \chi^{[2]}} \Z \to \Z/2^{[3]}$ the universal  4-bundle over $K(\E_{b \times \bb} \int_{\hat \chi^{[2]}} \Z,1)$
is well-defined, denote this bundle by 
$\hat \tau_{b \times \bb}$. The bundle 
\begin{eqnarray}\label{tauN}
\hat \eta_{b \times \bb\circ}^{\ast}(\hat \tau_{b \times \bb})
\end{eqnarray}
over $\hat N_{b \times \bb\circ}$ is well-defined.  

Denote by
\begin{eqnarray}\label{NNbb}
\widehat {NN}_{\circ} \subset \hat N_{b \times \bb \circ}
\end{eqnarray}
the 3-dimensional subpolyhedron, generally speaking, with boundary, as a homology Euler class
of the Whitney sum of 
$\frac{n-2k-3}{4}$ copies of the bundle $(\ref{tauN})$. The condition
 $(\ref{usll4})$ is equivalent to the following: the homology class 
\begin{eqnarray}\label{usll5}
(p_{\E_{d} \times \Z,\E_d} \circ \tilde \eta_{b \times \bb\circ})_{\ast}([C\widetilde {NN}_{\circ}]) \in  H_{3}(K(\E_d,1);\Z)
\end{eqnarray}
is even.

Consider the mapping
$\widehat {NN}_{\circ} \to K(\E_{b \times \bb} \int_{\chi^{[2]}} \Z,1) \to K(\Z,1)$.
Without loss of the generality, the inverse image by this mapping of the marked point of
$S^1 = K(\Z,1)$ is a closed 2-dimensional subpolyhedron, denoted by  
\begin{eqnarray}\label{LL}
\widehat {LL} \subset \widehat {NN}_{\circ}.
\end{eqnarray}
This polyhedron is ${\rm{PL}}$--homeomorphic to an oriented surface, which is equipped with a mapping
\begin{eqnarray}\label{fLL}
\hat f: \widehat {LL} \longrightarrow K(\E_{b \times \bb},1).
\end{eqnarray}
Let us use the following isomorphism: 
$H_2(K(\E_{b \times \bb},1);\Z)=\Z/2$. 

Let us prove that there exists a closed oriented 3-manifold
$\widehat {NN}$, its submanifold as in the formula  $(\ref{LL})$ and a mapping  
\begin{eqnarray}\label{hatF}
\hat F: \widehat {NN} \to K(\E_{b \times \bb} \int_{\chi^{[2]}} \Z,1), 
\end{eqnarray}
for which the following two conditions are satisfied:

--1. The image of the fundamental class by the mapping 
$(\ref{fLL})$ determines the generator of the group  $H_2(K(\E_{b \times \bb},1);\Z)$.

--2. The image of the fundamental class by the mapping 
$$\tilde F: \widetilde {NN} \to K(\E_d \times \Z,1) \to K(\E_d,1) = K(\Z/4,1) $$
is an even (or the trivial) element in the group 
$H_3(\E_d;\Z)$. 

Let us consider 2-torus $\widehat {LL}$, which is the the 2-skeleton of the standard cell decomposition of the space
$(\RP^{\infty} \times \RP^{\infty})/T_{\i} \supset (\RP^1 \times \RP^1)/T_{\i} = \widetilde {LL}$, where 
$T_{\i}:\RP^{\infty} \times \RP^{\infty} \to \RP^{\infty} \times \RP^{\infty}$ is the diagonal involution, 
which is defined by the standard involution $\i: \RP^{\infty} \to \RP^{\infty}$.  We may visualized the space $K(\E_d,1)$  as the space
$(\RP^{\infty} \times \RP^{\infty})/T_{\i} \setminus diag(\RP^{\infty})$.
By this construction the involution
$\hat \chi^{[2]}: K(\E_d,1) \to K(\E_d,1)$, which corresponds to the automorphism  $(\ref{hatchiE})$
is defined by the formula:  $x \times y \mapsto y \times x$. 

Define the (orientation preserving) involution $\hat{\chi}: \widehat {LL} \to \widehat {LL}$, which permutes the factors and reverses the diagonal. 
Define the mapping $\hat f: \widehat {LL} \to  K(\E_{b \times \bb},1)$ $(\ref{fLL})$, which transforms the diagonal 
generator 
$\i \in H_1(\widehat {LL};\Z)$ to the element $ab \in \H_{b \times \bb}$ (this element is represented by the sum of the diagonal loop with the generic loop of the first factor). Obviously, the mapping  $\hat f$ 
commutes up to homotopies with the involutions $\hat \chi$, $\hat \chi^{[2]}$ in the source and target spaces of the mapping
$\hat f$. Let us call the considered property Gluing Condition.

 Let us define the manifold
$\widehat {NN}$ as an oriented 3-manifold by the cylinder of the involution $\hat \chi: \widehat {LL} \to 
\widehat {LL}$.
The mapping  $(\ref{hatF})$ is well-defined by a fibered family over $S^1$ of mappings of
2-tori in the space
$K(\E_{b \times \bb},1)$ (the source and the target space of $(\ref{hatF})$ is the total spaces of fibrations
over $S^1$). By Gluing Condition the mapping $(\ref{hatF})$ is well-defined. This mapping satisfies Condition 1.

Let us check Condition 2. Consider the following composition:
\begin{eqnarray}\label{NNZ2}
p_{\E_d,\Z/2} \circ \tilde F: \widetilde {NN} \to K(\E_{d} \times \Z,1) \to K(\E_d,1) \to K(\Z/2,1),
\end{eqnarray}
where the mapping
$p_{\E_d,\Z/2}: K(\E_d,1) \to K(\Z/2,1)$ is induced by the epimorphism   $\E_d \to \Z/2$ with the kernel $\I_d \subset \E_d$. It is well-known, that the cellular mapping
 $p_{\E_d,\Z/2}$ transforms the standard 3-skeleton $S^3/\i \subset K(\E_d,1)$ into the standard 3-skeleton 
$\RP^3 \subset K(\Z/2,1)$ with degree 2.

Assuming Condition 2 is not satisfied and the mapping 
 $(\ref{hatF})$ determines the generic homology class, then the mapping
$(\ref{NNZ2})$ is not homotopic to zero. Assume that the mapping $(\ref{NNZ2})$ is cellular. 
Then the image of this mapping coincides with the standard 3-skeleton  $\RP^3 \subset K(\Z/2,1)$
and the degree of the mapping $(\ref{NNZ2})$ is equal to $2$ modulo 4. 
 
The mapping $(\ref{NNZ2})$ is a 2-sheeted covering over the mapping 
\begin{eqnarray}\label{hatNNZ2}
\widehat {NN} \to K(\E_{b \times \bb} \int _{\chi^{[2]}} \Z,1) \to K(\Z/2 \times \Z,1) \to K(\Z/2,1).
\end{eqnarray}
By the construction, the mapping 
 $(\ref{hatNNZ2})$ is homotopic to a mapping into the standard 2-skeleton $\RP^2 \subset K(\Z/2,1)$. 
This implies that image of the fundamental class by the mapping  
$(\ref{hatNNZ2})$, and by the mapping 
 $(\ref{NNZ2})$ is the trivial homology class. This prove that the degree of the mapping
 $(\ref{NNZ2})$ is equal to $0$ modulo 4. The mapping $\hat F$ satisfies Condition 2.

To prove Condition 
$(\ref{usll5})$ we may assume that the image of the fundamental class by the mapping $(\ref{fLL})$
is the trivial homology class. Therefore it is sufficiently to prove Condition $(\ref{usll5})$, assuming, that the surface $\hat LL$ is empty. In this case the mapping  $\hat \eta_{b \times \bb\circ}$ admits a reduction into the
subspace
 $K(\E_{b \times \bb},1) \subset K(\E_{b \times \bb}\int_{\chi^{[2]}} \Z,1)$.
Condition $(\ref{usll5})$ is reformulated analogously to Condition  $(\ref{usll2})$, which was proved above.
Condition 2 from [Lemma 27, A1] is proved. Lemma 
 $\ref{osnlemma1}$A  is proved.

\section{Proof of Lemma $\ref{7}$. Sketches of proofs of Lemma $\ref{osnlemma2}$, Proposition $28$ ${\rm{[A2]}}$, Proposition  $31$ ${\rm{[A2]}}$
and Lemma  $35$ ${\rm{[A2]}}$}

To prove Lemma
$\ref{7}$ is sufficiently to repeat a part of Lemma $\ref{osnlemma1}$ B., which is related with a subpolyhedron $RK_{b\times \bb\circ}$ in the polyhedron of the self-intersection. 
Lemmas
$30$, $32$ from ${\rm{[A2]}}$ are proved analogously to $\ref{osnlemma1}$. Lemma  35 from ${\rm{[A2]}}$
is proved analogously to Lemma
$\ref{osnlemma2}$. A detailed proof of the lemmas  requires to make the paper greater.

\subsubsection*{A sketch of the proof of Lemma  $\ref{osnlemma2}$}

The proof is analogous to the proof of the main result of the paper
 $\cite{Akh1}$. Let us consider an auxiallary mapping
$p_1: S^{n-2k+n_{\sigma-1}+1}/\i \to J_1$,
given by the
formula  $(\ref{p_1})$, define by $C_{p_1}$ the cylinder of this
mapping. The projections  $\pi_I: C_{p_1} \to [0,1]$, $\pi_J:
C_{p_1} \to J_1$ are well defined, denote the Cartesian product of
this mappings by $F_1: C_{p_1} \to J_1 \times [0,1]$.

Аналогично рассмотрим отображение $\tilde p_1: S^{n-2k}/\i \to
J_1$, определенное по формуле $(\ref{tildep_1})$ и обозначим через
$C_{\tilde p_1}$ цилиндр этого отображения. Определены отображения
проекций $\tilde \pi_I: C_{\tilde p_1} \to [0,1]$, $\tilde \pi_J:
C_{\tilde p_1} \to J_1$ и декартово произведение этих отображений,
которое обозначим через $\tilde F_1:  C_{\tilde p_1} \to J_1
\times [0,1]$. Определено вложение $r_1: C_{\tilde p_1} \subset
C_{ p_1}$. Cледующие диаграммы коммутативны:

\begin{eqnarray}\label{A}
\begin{array}{ccc}
C_{\tilde p_1} & \longrightarrow & C_{p_1} \\
\downarrow \tilde \pi_I & & \swarrow \pi_I  \\
I & & \\
\end{array}
\end{eqnarray}

\begin{eqnarray}\label{B}
\begin{array}{ccc}
C_{\tilde p_1} & \longrightarrow & C_{p_1} \\
\downarrow \tilde \pi_J & & \swarrow \pi_J  \\
J_1 & & \\
\end{array}
\end{eqnarray}

Consider the inclusion  $I_J: J_1 \times [0,1] \subset \R^n \times
[0,1]$ and define the mapping  $I_J \circ \tilde F_1: C_{\tilde
p_1} \to \R^n \times [0,1]$, $I_J \circ F_1: C_{p_1} \to \R^n
\times [0,1]$. Consider the mapping $\tilde f_1: C_{\hat p_1} \to
\R^n \times [0,1]$ which was defined by a small generic alteration
of the mapping  $I_J \circ \tilde F_1$. The mapping $\tilde f_1$
will be taken to be coincided on the bottom of the cylinder $J_1
\subset C_{\tilde p_1}$ with the embedding  $I_J: J_1 \subset \R^n
\times \{0\}$. Moreover, the composition $p_{[0,1]} \circ \tilde
f_1 : C_{\tilde p_1} \to [0,1]$ to be coincided with  $\tilde
p_{I}$, where $p_{I}: \R^n \times [0,1] \to [0,1]$ is the
projection on the second factor. The mapping $f: C_{p_1} \to \R^n
\times [0,1]$ is also defined such that $\tilde f_1 = f_1 \circ
r_1$.

Denote by $\bar Q_1 \subset C_{p_1}$ the polyhedron of
self-intersection points of the mapping  $f_1$, defined as the
closure of the corresponded spaces by the formula:
$$ \bar Q_1 = Cl\{ x \in C_{p_1} : \exists y \in C_{p_1}, x \ne y,
f(x) = f(y) \}. $$ Because $n-4k=n_{\sigma}$, $\dim(\bar {\tilde
Q}_1)=n_{\sigma+1}+1$.

Denote by $\bar {\tilde Q}_1 \subset C_{\tilde p_1}$ the
polyhedron of self-intersection points of the mapping  $\tilde
f_1$, this polyhedron is defined as the closure of the
corresponded subspaces by the formula
$$ \bar {\tilde Q}_1 = Cl\{ x \in C_{\tilde p_1} : \exists y \in C_{\tilde p_1}, x \ne y, \tilde
f_1(x) = \tilde f_1(y) \}. $$ Because $n-4k=n_{\sigma}$, we get
$\dim(\bar Q_1)=n_{\sigma}+1$.

Consider the stratification $J_1^{[2]} \subset J_1^{[1]} \subset
J_1$ of the join. Denote by $\bar Q_{J_1}$ the intersection $\bar
Q \cap J_1$. Denote by  $\bar {\tilde Q}_{J_1}$ the intersection
$\bar {\tilde Q}_1 \cap J_1$. The polyhedron  $\bar Q_{J_1}$ has
the codimension $n_{\sigma+1}$. Because the codimension of
$J^{[2]}_1 \subset J_1$ is equal to $n_{\sigma+1}+1$, the
polyhedron  $\bar Q_{J_1} \subset J_1$ is outside a regular
neighborhood of the stratum $J^{[2]}_1$. The polyhedron  $\bar
{\tilde Q}_{J_1}$ has the codimension $n_{\sigma}$. Because the
codimension of  $J^{[1]}_1 \subset J_1$ is equal to
$n_{\sigma}+1$, the polyhedron $\bar {\tilde Q}_{J_1} \subset J_1$
is outside a regular neighborhood of the stratum  $J^{[1]}_1$.
Define the polyhedron $\bar {\tilde Q}_{J_1}(\varepsilon)$ as the
set of points from   $\bar {\tilde Q}_{J_1}$ which are mapped with
respect to the projection $\tilde \pi_I$ into a small positive
$\varepsilon \in I$.

Define the involution  $T_{\bar {\tilde Q}}: \bar {\tilde Q} \to
\bar {\tilde Q}$ which permutes points of self-intersection on the
canonical covering. The involution   $T_{\bar {\tilde Q}}$ keeps
the values of the mapping $\tilde \pi_I$. The polyhedron $\bar
{\tilde Q}_{J_1}(\varepsilon)$ is invariant with respect to the
involution   $T_{\bar {\tilde Q}}$. Denote by $T_{\bar {\tilde
Q}}(\varepsilon)$ the restriction of the considered involution on
the polyhedron $\bar {\tilde Q}_{J_1}(\varepsilon)$, this
restriction is a free involution.

Define the mapping $d_1: S^{n-2k}/\i \to \R^n \times
\{\varepsilon\} = \R^n$ as the restriction of the mapping  $\tilde
f_1$ on $S^{n-2k}/\i \times \{\varepsilon\}$. A quotient $\bar
{\tilde Q}_{J_1}(\varepsilon)/T_{\bar {\tilde Q}}(\varepsilon)$ is
a polyhedron of self-intersection points of the mapping  $d_1$.
Consider the polyhedron of self-intersection of the mapping
 $d_1$ and  its subpolyhedron  $N_1$.
By the construction, if the positive parameter
$\varepsilon$ is small enough,  the structured mapping $\zeta: N_1
\to K(\H,1)$ admits a reduction to a mapping into the subspace
 $K(\Q,1) \cup
K(\H_b,1) \subset K(\H,1)$, the considered reduction is well
defined as the composition of the mapping $t_1: N_1 \to RK_1$ with
the mapping $\phi_1: RK_1 \to K(\Q,1) \cup K(\H_b,1)$ (see the diagram $(\ref{RK1})$). 

Let us prove that the mapping
 $t_1$ satisfies the boundary conditions from diagram  $(\ref{118.21})$ in Lemma $\ref{lemma291}$. 
For $\ell \ge 8$ the number  $r_1$ of the factors of the join $J_1$, which is calculated by the formula 
$(\ref{r_1})$, is greater then  $n_{\sigma}$. Because $\dim(N_1)=n_{\sigma-1}-1$, the boundary of the polyhedron
$N_1$ contains no strata of a deep greater then  $\frac{r_1-1}{2}$. Therefore the coordinate system 
in each component  $N_1$ of the type $\E_b$ is agree with boundary conditions. Lemma  $\ref{osnlemma2}$ is proved.

\subsubsection*{$\H_{b \times \bb}$--structure of formal mappings with holonomic singularities}

Consider the polyhedron 
$X_{b \times \bb} \int_{\chi} S^1$, which is a skeleton of the Eilenberg-Mac Lane
space 
$K(\I_{b \times \bb} \int_{\chi^{[2]}} \Z),1)$, correspondingly to [Formula (181), A2]. 
Consider the mapping ${\rm{i}}_{J_X \int_{\chi} S^1} \circ \varphi_{X_{b \times \bb}}: 
X_{b \times \bb} \int_{\chi} S^1 \to D^{n-1} \times S^1 \subset \R^n$, where the mapping
$\varphi_{X_{b \times \bb}}$ is defined by the [Formula (186), A2], and the mapping (embedding) 
${\rm{i}}_{J_X \int_{\chi} S^1}$ is defined by the [Formula (190), A2]. We shall consider this mapping 
as a mapping with a holonomic singularity in the sense of [Definition 
 9, A2]. Denote this formal mapping by  $(d_{\int,0}, d_{\int,0}^{(2)})$. Let us restrict this formal mapping  $(d_{\int,0},d_{\int,0}^{(2)})$ on the subpolyhedron
$X_{b \times \bb} \subset X_{b \times \bb} \int_{\chi} S^1$, and denote this restriction by  
$(d_0,d_0^{(2)})$.

\begin{lemma}\label{lemmaE}
There exists a $C^0$--small  ${\rm{PL}}$-deformation of the formal holonomic pair of mappings  
$(d_{\int,0},d_{\int,0}^{(2)})$ to a pair of mappings $(d_{\int},d^{(2)}_{\int})$ 
with holonomic singularity, such that  
the
polyhedron $N_{\int\circ}$ of formal self-intersection of the mapping  $(d_{\int},d^{(2)}_{\int})$ is decomposed into the union of two subpolyhedra:
\begin{eqnarray}\label{bb[3]}
N_{\int \circ} = N_{\int, \H_{b \times \bb}} \cup N_{\int, [3] \circ},
\end{eqnarray}
where $N_{\int, \H_{b \times \bb}}$ is closed. 

The restriction of the structure mapping
 $\zeta_{\circ}$ on the subpolyhedron  $N_{\int, \H_{b \times \bb}}$  admits a reduction, which is given by the mapping
 $\zeta_{b \times \bb}: N_{\int, \H_{b \times \bb}} \to K(\H_{b \times \bb} \int _{\chi^{[3]}} \Z,1)$.



The polyhedron
$N_{\int \circ}$ contains a subpolyhedron  $N_{\circ} \subset N_{\int \circ}$, 
which is decomposes into two components:
$$N_{\circ} =N_{\H_{b \times \bb}} \cup N_{[3]\circ},$$
where the components are defined as the corresponding components in the formula 
$(\ref{bb[3]})$.
The restriction of the structured map $\zeta_{\circ}$ on the subpolyhedron  $N_{[3] \circ}$ 
admits a reduction, which is given by the mapping   
 $$ \zeta_{b \times \bb \times \Z/2\circ}: N_{[3]\circ} \to K((\I_{b \times \bb} \times \Z/2) \int_{\chi} \Z,1),$$
and which is satisfies the boundary condition, given by a mapping into the subspace $K(\I_{b \times \bb} \times \Z/2,1)$.
 (In this formula the extension of the group $\I_{b \times \bb} \times \Z/2$ (and analogous extensions below) 
 are corresponding to the inclusion  $X_{b \times \bb} \subset X_{b \times \bb} \int_{\chi} S^1$.)
 
 The mapping  $\zeta_{b \times \bb \times \Z/2\circ}$ is a compressed by the canonical 2-sheeted covering 
$N_{[3]\circ} \to \hat N_{[3]\circ}$, and is a 2-sheeted covering mapping over the mapping 
  $$\hat \zeta_{b \times \bb \times \Z/2 \circ}: \hat N_{[3] \circ} \to K((\H_{b \times \bb} \times \Z/2) \int_{\hat \chi} \Z,1),$$
 which is satisfies the boundary condition, given by a mapping into the subspace
$K((\H_{b \times \bb} \times \Z/2)\int_{\hat \chi} \Z,1)$.  In the previous formula the automorphism  (involution)
$\hat \chi:  \H_{b \times \bb} \times \Z/2 \to  \H_{b \times \bb} \times \Z/2$ is the identity on the subgroup
$\H_{b \times \bb} \subset \H_{b \times \bb} \times \Z/2$, and is mapped the generator 
$t \in \Z/2$ into the element  $t_d t$, where $t_d$ is the generator of the subgroup $\I_d \subset \H_{b \times \bb}$.
Define the automorphism $\chi:  \I_{b \times \bb} \times \Z/2 \to  \I_{b \times \bb} \times \Z/2$
by the restriction of 
$\hat \chi$ on the subgroup.  
\end{lemma}

Let us formulated and proof a lemma, which is required to check [Formula
 (211), A2]. For an arbitrary pair of integers
$(s_1,s_2)$, $s_1 = 1
\pmod{2}$, $s_2 = 1 \pmod{2}$,  $s=s_1 + s_2 = n - \frac{n-m_{\sigma}}{8}$,
consider the homology class [(210), A2].  This homology class is defined as the image of the fundamental class
of the manifold
 $X(s_1,s_2)$, which is naturally embedded into
$X_{b \times \bb}$.

Denote the restriction of $(d,d^{(2)})$ on
$X(s_1,s_2)$ by $(d(s_1,s_2),d^{(2)}(s_1,s_2))$. Consider a polyhedron of the formal self-intersection 
of the mapping
 $(d(s_1,s_2),d^{(2)}(s_1,s_2))$, which is represented by a disjoin union of the two subpolyhedra. 
The canonical covering over the first polyhedron is a closed subpolyhedron into  
$\bar {N}_{\H_{b \times\bb}}$, the canonical covering over the second polyhedron 
is the closure of an open subpolyhedron in 
$C\bar {N}_{[3]\circ}$, denote this closure by $C\overline{NX}(s_1,s_2)$. 

Denote the fundamental class of the polyhedron
 $C\overline{NX}(s_1,s_2)$ by
$[C\overline{NX}(s_1,s_2)] \in H_{n-\frac{n-m_{\sigma}}{4}}(K(\I_{b \times \bb},1))$
 (we have used the isomorphism [(42),A2]). 

Let us prove the formulas
[(211),A2] analogously to Lemma
$\ref{osnlemma1}$A, in which the formula  $(\ref{usll})$ is proved (Condition 1 from [Lemma 26, A1]). 

\begin{proposition}\label{prop220}
An arbitrary homology class
$[C\overline{NX}(s_1,s_2)]$ is trivial.
\end{proposition}

\subsubsection*{Proof of Proposition $\ref{prop220}$}
Denote by 
 $NX(s_1,s_2)_{\circ}$ an open polyhedron, which is the base of 2-sheeted covering
 space  $\overline{NX}(s_1,s_2)_{\circ}$.  The polyhedron
$NX(s_1,s_2)_{\circ}$ is equipped with the structure mapping
$$ \zeta(s_1,s_2)_{\circ}: NX(s_1,s_2)_{\circ} \to K((\I_{b \times \bb} \times \Z/2) \int_{\hat \chi} \Z,1), $$
and the regular neighborhood of the boundary is mapped by the considered structure mapping
into the subspace 
\begin{eqnarray}\label{subspaceI}
K(\I_{b \times \bb} \times \Z/2,1) \subset K((\I_{b \times \bb} \times \Z/2) \int_{\hat \chi} \Z,1).
\end{eqnarray}

The manifold
$X(s_1,s_2)$ is a 2-sheeted covering over the manifold $\hat X(s_1,s_2)$.
Therefore, an open polyhedron, which is a base of the 2-sheeted covering  with the covering space 
$NX(s_1,s_2)_{\circ}$ is well-defined. Let us denote this polyhedron by  $\widehat{NX}(s_1,s_2)_{\circ}$. 
The polyhedron  $\widehat{NX}(s_1,s_2)_{\circ}$ is equipped with a structure mapping
$$ \hat{\zeta}_{\circ}: \widehat{NX}(s_1,s_2)_{\circ} \to K((\H_{b \times \bb} \times \Z/2) \int_{\hat \chi} \Z,1), $$
a regular neighborhood of the boundary is mapped by this mapping into the subspace
\begin{eqnarray}\label{subspace}
K(\H_{b \times \bb} \times \Z/2,1) \subset K((\H_{b \times \bb} \times \Z/2) \int_{\hat \chi} \Z,1).
\end{eqnarray}

Assume, that the image of the structure mapping
$\zeta(s_1,s_2)_{\circ}$ is inside the subspace  $(\ref{subspaceI})$.
Then the statement of the lemma is evident, because 
$C\overline{NX}(s_1,s_2)$ is a composition with a 2-sheeted covering over 
$CNX(s_1,s_2)$ (comp. with the initial step of the proof of Lemma $\ref{osnlemma1}$A). 

Let us consider a general case. We shall use the polyhedron 
$\widehat{NX}(s_1,s_2)_{\circ}$. The universal bundle over the space 
$K((\H_{b \times \bb} \times \Z/2) \int_{\hat \chi} \Z,1)$ is a 8-dimensional bundle.
It is sufficiently to prove the formula for the cycle, which is defined as the intersection 
of the considered  fundamental class with the Euler class of the pull-back of a suitable Whitney sum of the universal bundle.
Denote the Euler class of the universal bundle by
$\hat \tau_{b \times \bb \times \Z/2,\int}$. 

For an arbitrary pair of the positive integers
$(p_1,p_2)$, $p_1 = 1
\pmod{2}$, $p_2 = 1 \pmod{2}$,  $p=p_1 + p_2 = \frac{n + 6}{2}$,
define the submanifold
 $XX(p_1,p_2)$ of the dimension $p$, $XX=\RP^{p_1} \times \RP^{p_2}$.

Define the embedding
 $XX(p_1,p_2) \subset X(s_1,s_2)$, as the Cartesian product of the coordinate embeddings 
$\RP^{p_1} \subset \RP^{s_1}$, $\RP^{p_2} \subset \RP^{s_2}$, which satisfies the restriction $s_1-p_1 = s_2-p_2 = \frac{s-p}{2}$. Define the formal mapping  $(dd(p_1,p_2),dd^{(2)}(p_1,p_2))$ as the restriction of the formal mapping $(d(s_1,s_2),d^{(2)}(s_1,s_2))$ to the submanifold
$XX(p_1,p_2)$. Denote by $NXX(p_1,p_2)_{\circ}$ an open polyhedron of the formal self-intersection 
of the mapping 
 $(dd(p_1,p_2),dd(p_1,p_2)^{(2)})$. The following 6-dimensional subpoluhedron
$$NXX(p_1,p_2)_{\circ} \subset NX(s_1,s_2)_{\circ}$$
is well-defined, the fundamental class of this subpolyhedron is realized the homology Euler class
of the bundle
$\zeta_{\circ}^{\ast}(\tau_{b \times \bb \times \Z/2, \int}^{\frac{s-p}{8}})$.

Let us prove that the homology class
\begin{eqnarray}\label{CNXX}
[C\overline{NXX}(p_1,p_2)] \in H_{6}(\I_{b \times \bb},1)
\end{eqnarray}
 is trivial. We shall distinguishes the exceptional case,  when  
$p_1=1$, or $p_2=1$.  Consider non-exceptional case in which  $p_1 \ge 3$, $p_2 \ge 3$. 
Let us prove that the homology class
 $(\ref{CNXX})$ is trivial. 

The lens manifold
$(\RP^{p_1} \times \RP^{p_2})/\i_{diag}$
is immersible into $\R^n$.  Therefore the homology class of the boundary singularities of the polyhedron 
$\partial(\widehat{NXX}(p_1,p_2)_{\circ})$ in the group  $H_5(\H_{b \times \bb} \times \Z/2,1)$ 
is trivial. Let us omit below the marks 
 $\circ$ and $C$ in denotations.

Let us consider the 5-dimensional fundamental class
$[\hat{p}^{-1}(pt)] \in H_5(\H_{b \times \bb} \times \Z/2,1)$ of the closed subpolyhedron  $\hat{p}^{-1}(pt)$, where
$\hat{p}: \widehat{NXX}(p_1,p_2) \to S^1$ is the projection, which is induced by the projection 
$p_{\H_{b \times \bb} \times \Z/2,\int}$ of the universal space. 

Assume that the homology class
 $[\hat p^{-1}(pt)]$  is trivial. Then, without loss of a generality, we may assume that the manifold
   $\hat p^{-1}(pt)$ is empty and the proof is reduced to the previous.  

Assume that the homology class
$[\hat p^{-1}(pt)]$ is non-trivial. Let us prove that the homology class
$[\hat p^{-1}(pt)]^!$ is realized for a suitable mapping of a closed 6-dimensional manifold
$A$,
$\zeta_A: A \to (\H_{b \times \bb} \times \Z/2) \int_{\hat \chi} \Z,1)$, for which the homology class,
defined analogously to  $(\ref{CNXX})$, is trivial. 

Let us decompose the fundamental class
$[\hat p^{-1}(pt)]$  over the base of the group 
$Im(H_5(K(\H_{b \times \bb} \times \Z/2,1);\Z) \to H_5(K(\H_{b \times \bb} \times \Z/2,1))$. 
Consider the following epimorphisms:
$$\pi_b: \H_{b \times \bb} \times \Z/2 \to \H_{b} \times \Z/2,$$
$$\pi_{\bb}: \H_{b \times \bb} \times \Z/2 \to \H_{\bb} \times \Z/2.$$

Assume that the image of the homology class
 $[\hat p^{-1}(pt)]$  in the group
$H_5(K(\H_b \times \Z/2,1) \times K(\H_{\bb} \times \Z/2,1))$ by the homomorphism  $(\pi_b \times \pi_{\bb})_{\ast}$
is represented by the tensor product of a homology class of
$H_2(K(\H_b \times \Z/2,1)$ to  a homology class of
$H_3(K(\H_{\bb} \times \Z/2,1))$. The proof in the last cases is evident (or is is given after $b$ is replaced by $\bb$.)

The condition 
 $\hat \chi_{\ast}([\hat p^{-1}(pt)])= [\hat p^{-1}(pt)]$ is satisfied, because the boundary conditions on 
$\widehat{NXX}(p_1,p_2)_{\circ}$ determines the trivial homology class. Therefore, after the expansion of the
element 
 $\pi_{\bb,\ast}([\hat p^{-1}(pt)])$ over the standard base the generator
of the factor $H_{3}(K(\Z/2,1);\Z)$ is not involved and  $\pi_{\bb,\ast}([p^{-1}(pt)])$
is expressed by the generator of
$H_3(K(\H_b,1))$.

Analogous to the construction 
 $(\ref{hatF})$, without loss of a generality, we may assume that the homology class  $(\ref{CNXX})$ 
 is trivial. Therefore, without loss of a generality, we may assume, that   $p^{-1}(pt) = \emptyset$, 
 and we may repeat the previous proof as in the case, when the image of the structure mapping is inside the subspace 
  $(\ref{subspace})$. 

Is sufficiently to prove that in the exceptional case the homology class 
$(\ref{CNXX})$ is trivial. Let us  decomposes the homology class $(\ref{CNXX})$ over the standard base of the group 
 $H_{6}(\I_{b \times \bb},1)$. The generators of the group are 
$t_{3,b}t_{3,\bb}$, $t_bt_{5,\bb}$, $t_{5,b}t_{\bb}$. In the exceptional case, evidently, that the generator
$t_{3,b}t_{3,\bb}$ is not involved. To prove that the last generators
$t_bt_{5,\bb}$, $t_{5,b}t_{\bb}$ are not involved, let us intersect the 6-dimensional polyhedron  $NXX(p_1,p_2)_{\circ}$ with 4-dimensional Euler class of the universal bundle, which is the bull-back by $\pi_b$, 
or by $\pi_{\bb}$, correspondingly to the generators 
 $t_bt_{5,\bb}$, $t_{5,b}t_{\bb}$. The proof is analogous to the previous proof, this proof is more simple, because the Euler class is represented by a 2-dimensional subpolyhedron in  $NXX(p_1,p_2)_{\circ}$. 
 Is sufficiently to consider the only generators of 
 $H_1(K(\H_{b \times \bb} \times \Z/2,1);\Z) = \H_{b \times \bb} \times \Z/2$.
Lemma $\ref{prop220}$ is proved. 
\[ \]

 \[  \]
 Troitsk, IZMIRAN, 142190
\[  \]
 pmakhmet@izmiran.ru
\[  \]

\end{document}